\documentclass[reqno]{amsart}

\usepackage{amsmath,amssymb, amsopn, amsthm}
\usepackage{enumitem, framed, xcolor}
\usepackage{changes}
\usepackage{tikz}

\usetikzlibrary{arrows.meta,calc,positioning}

\definecolor{annred}{RGB}{238,48,71}
\definecolor{anngreen}{RGB}{40,205,92}
\definecolor{annblue}{RGB}{30,145,245}

\tikzset{
  fin/.style={draw=black, dash pattern=on 5pt off 4pt, line width=.75pt},
  nonfin/.style={draw=black, line width=2.15pt},
  unknown/.style={draw=black, line width=.55pt},
  annfinred/.style={draw=annred, dash pattern=on 4pt off 4pt, line width=1.0pt},
  annfingreen/.style={draw=anngreen, dash pattern=on 4pt off 4pt, line width=1.0pt},
  annfinblue/.style={draw=annblue, dash pattern=on 4pt off 4pt, line width=1.0pt},
  annnot/.style={draw=annblue, line width=2.0pt},
  op/.style={font=\normalsize, inner sep=1pt},
  edgelab/.style={font=\small, inner sep=1pt, fill=white},
  smalllab/.style={font=\small, inner sep=1pt},
  anntext/.style={font=\small\sffamily, align=left}
}

\theoremstyle{definition} 
\newtheorem{theorem}{Theorem}[section]
\newtheorem{lemma}[theorem]{Lemma}

\newtheorem{remark}[theorem]{Remark}

\newtheorem{problem}[theorem]{Problem}

\newtheorem{definition}[theorem]{Definition}

\renewenvironment{proof}[1][\unskip]{%
\par
\noindent
\textbf{Proof #1.}
\noindent}
{\hfill$\blacksquare$

\bigskip}

\let\<=\langle
\let\>=\rangle

\def\[#1\]{\begin{align*}#1\end{align*}}

\let\models=\vDash

\renewcommand{\setminus}{\smallsetminus}

\def\defeq{=}

\def\twoheadrightarrowtail{%
  \mathrel{\ooalign{%
    $\rightarrowtail$\cr
    \hidewidth$\twoheadrightarrow$\hidewidth\cr
  }}%
}
\def\izomorf{\twoheadrightarrowtail}

\DeclareMathOperator{\ran}{ran}
\DeclareMathOperator{\dom}{dom}

\DeclareMathOperator{\id}{id}

\def\gA{\mathfrak{A}}

\def\gB{\mathfrak{B}}

\makeatletter
\newcommand\RedeclareMathOperator{%
  \@ifstar{\def\rmo@s{m}\rmo@redeclare}{\def\rmo@s{o}\rmo@redeclare}%
}
\newcommand\rmo@redeclare[2]{%
  \begingroup \escapechar\m@ne\xdef\@gtempa{{\string#1}}\endgroup
  \expandafter\@ifundefined\@gtempa
     {\@latex@error{\noexpand#1undefined}\@ehc}%
     \relax
  \expandafter\rmo@declmathop\rmo@s{#1}{#2}}
\newcommand\rmo@declmathop[3]{%
  \DeclareRobustCommand{#2}{\qopname\newmcodes@#1{#3}}%
}
\@onlypreamble\RedeclareMathOperator
\makeatother

\DeclareMathOperator{\Sb}{Sb}
\DeclareMathOperator{\Rd}{Rd}
\DeclareMathOperator{\SRd}{SRd}

\DeclareMathOperator{\C}{C}
\RedeclareMathOperator{\S}{S}
\DeclareMathOperator{\D}{D}
\RedeclareMathOperator{\P}{P}
\DeclareMathOperator{\rep}{rep}

\def\FSC{\mathsf{SCA}}
\def\RCA{\mathsf{RCA}}
\def\BA{\mathsf{BA}}
\def\PA{\mathsf{PA}}
\def\RA{\mathsf{RA}}

\def\FPA{\mathsf{FPA}}
\def\FPEA{\mathsf{FPEA}}
\def\FRPA{\mathsf{RFPA}}
\def\FRPEA{\mathsf{RFPEA}}
\def\RPA{\mathsf{RPA}}
\def\RDF{\mathsf{Rdf}}
\def\RPEA{\mathsf{RPEA}}

\def\QPEA{\mathsf{QPEA}}
\def\BAO{\mathsf{BAO}}
\def\CA{\mathsf{CA}}
\def\PEA{\mathsf{PEA}}
\def\Df{\mathsf{Df}}
\def\QA{\mathsf{QA}}

\def\Str{\mathsf{Str}}

\def\c{\mathsf{c}}
\def\s{\mathsf{s}}
\def\p{\mathsf{p}}
\def\d{\mathsf{d}}

\def\cS{\mathcal{S}}
\def\Cm{\mathfrak{Cm}}
\def\Uf{\mathfrak{Uf}}
\def\uf{{U\!f}}
\def\Em{\mathfrak{Em}}

\let\phi=\varphi

\def\Bl{\mathfrak{Bl}}
\def\Sg{\mathfrak{Sg}}
\def\RfAt{RfAt}

\begin{document}

\title[Permutations, substitutions and finite axiomatizability]{Permutations, substitutions and finite axiomatizability}

\author[Andr\'eka]{Hajnal Andr\'eka}

\address{%
HUN-REN Alfr\'ed R\'enyi Institute of Mathematics \\
Re\'altanoda utca 13-15 \\
1053 Budapest \\
Hungary
}

\email{andreka.hajnal@renyi.hu}


\author[Gyenis]{Zal\'an Gyenis}
\address{%
	Jagiellonian Univerisity \\
	Grodzka 52 \\
	31-007 Krak\'ow \\
	Poland}

\email{zalan.gyenis@uj.edu.pl}

\author[N\'emeti]{Istv\'an N\'emeti}
\address{%
HUN-REN Alfr\'ed R\'enyi Institute of Mathematics \\
Re\'altanoda utca 13-15 \\
1053 Budapest \\
Hungary
}
\email{nemeti.istvan@renyi.hu}




\date{\today}

\begin{abstract}
	Algebras of relations form an algebraic framework for the study of logical systems, extending the correspondence between Boolean algebras and propositional logic. Tarski's representable cylindric algebras $\RCA_{\alpha}$, and Halmos' representable polyadic algebras $\RPA_{\alpha}$ both provide algebraic counterparts to first-order logic. 
	In this paper, we show that the usual finite set of polyadic axioms axiomatize $\RPA_{\alpha}$ over $\RDF_{\alpha}$, the diagonal-free subreducts of elements in $\RCA_{\alpha}$. In short: $\RPA_{\alpha} = \PA_{\alpha} + \RDF_{\alpha}$.
\end{abstract}

\dedicatory{Dedicated to Robert Goldblatt}

\maketitle

\section{What was and what is to come}\label{sec:1}

Algebras of relations have played an important role in the development of logic and its surrounding fields ever since the pioneering works of De Morgan, Peirce, and Schröder initiated the area. Boolean set algebras (that is, algebras of unary relations) are the algebraic counterparts of propositional logic: the equations valid in Boolean set algebras correspond to propositional tautologies. Alfred Tarski and Paul Halmos introduced the classes $\RCA_\alpha, \RDF_\alpha$ and $\RPEA_\alpha, \RPA_\alpha$ of algebras isomorphic to cylindric and polyadic algebras of $\alpha$-ary relations, respectively, as candidates for playing the same role in first-order logic (with and without equality) using $\alpha$ many variables as Boolean set algebras play in propositional logic. 
The class $\BA$ of Boolean algebras, the class of all algebras isomorphic to Boolean set algebras, is axiomatized by finitely many equations.
Similarly, all of the above-mentioned classes of cylindric and polyadic algebras are axiomatizable by equations if $\alpha$ is finite. However, unlike $\BA$, none of them can be axiomatized by finitely many equations if $\alpha\ge 3$.%
\footnote{For $\RCA_\alpha$, non-finite axiomatizability is proved by J.\ Donald Monk \cite{M69}, and for the rest by James Johnson \cite{J69}.
In contrast, for $\alpha<3$, all of them have finite axiom systems, see \cite[5.1.77, 3.2.65, 5.4.33, p.242]{HMT}.}\\

From now on, throughout the paper, $\alpha\ge 3$ is a finite ordinal. $T_{\alpha}$, or just  $T$, denotes the  set of all transformations $\sigma:\alpha\to\alpha$.\\

Cylindric and polyadic algebras are Boolean algebras with some extra operations. The similarity type of $\RPEA_\alpha$ 
has the following extra Boolean operations: the unary cylindrification operations  $\c_{(\Gamma)}$ and substitution operations $\s_{\sigma}$ for all subsets $\Gamma\subseteq\alpha$ and transformations $\sigma\in T$ together with the diagonal constants $\d_{ij}$ for all $i,j<\alpha$.
The other classes are obtained by discarding some of these operations and taking subalgebras: $\RPA_{\alpha}$ is the class of subalgebras of reducts of elements of $\RPEA_{\alpha}$ when discarding the diagonal constants $\d_{ij}, i,j<\alpha$; $\RCA_{\alpha}$ is the class of subalgebras of reducts of elements of $\RPEA_{\alpha}$ when discarding the substitution operations $\s_{\sigma},\ \sigma\in T$; and $\RDF_{\alpha}$ is the class of subalgebras of elements of $\RPEA_{\alpha}$ when we keep only the cylindrification operations. Hence,  $\RDF_{\alpha}$ is the class of all subalgebras of elements of $\RPA_{\alpha}$ with the substitution operations discarded.  See Figure~\ref{fig}.

Which of the extra Boolean operations introduce infinitely many new axioms, thereby causing non-finite axiomatizability? The cylindrifications clearly bring in infinitely many new axioms, since $\RDF_{\alpha}$ is non-finitely axiomatizable by \cite{J69} while $\BA$ is finitely axiomatizable. 
A natural question arises whether only the  cylindrification operations cause non-finite axiomatizability, or do the diagonal constants and the substitution operations contribute, too. 
Whether the diagonal constants or the substitution operations bring in infinitely many new axioms is asked in \cite{J69}. Namely, Problems 1,2 in \cite{J69} ask whether $\RCA_{\alpha}$ is finitely axiomatizable over $\RDF_{\alpha}$, and whether $\RPEA_{\alpha}$ is finitely axiomatizable over $\RPA_{\alpha}$ or over $\RCA_{\alpha}$.%
\footnote{If $K$ is a class of structures of similar types and $K'$ is a class of structures whose type is a reduct of the type of $K$, we say that $K$ is finitely axiomatizable over $K'$ if there is a finite set of axioms $\Gamma$ such that $\gA \in K$ iff $\gA$ is a model of $\Gamma$ and the reduct of $\gA$ to the type of $K'$ is a member of $K'$.}
All these questions have been answered negatively
in the meantime\footnote{In \cite[Thm.4]{AAPALI97ec}, \cite[Thm.6]{AAPALI97ec}, and \cite{ANTnew}, respectively.}, providing the answer to the question that all three kinds of extra-Boolean operations contribute to non-finite axiomatizability.
\begin{figure}[h]\begin{center}
		\begin{tikzpicture}[
			>=stealth,
			every node/.style={font=\normalsize},
			]
			
			\node (A) at (0,2) {$\mathsf{RPEA}_\alpha$};
			\node (B) at (6,2) {$\mathsf{RCA}_\alpha$};
			\node (C) at (0,0) {$\mathsf{RPA}_\alpha$};
			\node (D) at (6,0) {$\mathsf{Rdf}_\alpha$};
			
			\draw[->, ultra thick] (A) -- node[below] {$-\s_{\sigma}$} (B);
			
			\draw[->, ultra thick] (A) -- node[left] {$-\d_{ij}$} (C);
			
			\draw[->, ultra thick] (B) -- node[right] {$-\d_{ij}$} (D);
			
			\draw[->] (C) -- node[below] {$-\s_{\sigma}$} (D);
			
		\end{tikzpicture}
	\end{center}
	\caption{The three thick arrows represent relative non-finite axiomatizability. In the fourth line there is relative finite axiomatizability, this is proved in the present paper.}\label{fig}
\end{figure}

One natural question, however, was not asked in \cite{J69}. This is whether $\RPA_{\alpha}$ is finitely axiomatizable over $\RDF_{\alpha}$ or not. This question is answered in the present paper.
Surprisingly, we get an affirmative answer: $\RPA_{\alpha}$ is finitely axiomatizable over $\RDF_{\alpha}$. This is in sharp contrast with the situation when the diagonal constants are available: $\RPEA_{\alpha}$ can be axiomatized over $\RCA_{\alpha}$ with infinite and rather complex systems of equations only (see the main theorem of \cite{ANTnew}). 
Thus, we have to refine the picture: the substitution operations $\s_{\sigma}$ contribute to non-finite axiomatizability only when the diagonal constants are present in the similarity type, otherwise they do not contribute to non-finite axiomatizability.

Since all of $\RPEA_{\alpha},\RPA_{\alpha},\RCA_{\alpha},\RDF_{\alpha}$ are nonfinitely axiomatizable, finitely axiomatized ``approximating" classes $\PEA_{\alpha},\PA_{\alpha},\CA_{\alpha},\Df_{\alpha}$ were introduced for their study. In these varieties, the most important properties of the operations were gathered as axioms, and in such situations the question whether the finite axioms were well selected always is in the air. We  prove in this paper that the finitely many polyadic axioms introduced by Paul Halmos in the definition of $\PA_{\alpha}$ to describe the  behaviour of the substitution operations%
\footnote{See, for example, \cite[pp.27--28]{Halmos}, and \cite[Def.5.4.1]{HMT}.} are adequate in the sense that the equations (1)-(4) below together with the equational theory of $\RDF_{\alpha}$ axiomatize the equational theory of $\RPA_{\alpha}$: 
\begin{theorem}[$\s_{\sigma}$ is finitely axiomatizable over $\c_{(\Gamma)}$]\label{main1-thm} Assume that $3\le\alpha$ is finite.
Let  $\gA=\langle A,+,-,\c_{(\Gamma)},\s_{\sigma}\rangle_{\Gamma\subseteq\alpha, \sigma\in T}$. Then $\gA\in\RPA_{\alpha}$ if and only if
$\langle A,+,-,\c_{(\Gamma)}\rangle_{\Gamma\subseteq\alpha}\in\RDF_{\alpha}$ and the following equations hold in $\gA$ for all $\Gamma\subseteq\alpha$ and $\sigma,\tau\in T$.
\begin{align}
	&	\s_{\sigma}(x+y)=\s_{\sigma}(x)+\s_{\sigma}(y),\quad \s_{\sigma}(-x)=-\s_{\sigma}(x),\label{main1}\\
	&	\s_{\sigma\circ\tau}(x) = \s_{\sigma}\s_{\tau}(x),\quad \s_{Id}(x)=x,\label{main2}\\
	 &  \s_{\sigma}\c_{(\Gamma)}x=\s_{\tau}\c_{(\Gamma)}x\quad\mbox{if $\sigma,\tau$ differ only in $\Gamma$},\label{main3}\\
	&	\c_{(\Gamma)}\s_{\sigma}(x)=\s_{\sigma}\c_{(\sigma^{-1}\Gamma)}x,\quad\mbox{if $\sigma$ is one-one on $\sigma^{-1}\Gamma$}.	\label{main4}
\end{align}
\end{theorem}
\bigskip

\noindent
This is one of the main theorems of the present paper. We get Theorem~\ref{main1-thm} as a corollary to a more refined investigation, see Theorems \ref{main2-thm} and \ref{main3-thm} below. 
Let $[i\slash j]\in T$ and $[i,j]\in T$ denote the replacement and transposition operations, respectively, where $i,j<\alpha$. They are defined as, 
for every $k<\alpha$, 
\begin{align}
	[i/j](k) \defeq \begin{cases}
		j & \text{ if } k=i\\
		k & \text{ otherwise}
	\end{cases} \qquad
	\text{ and }\qquad
	[i,j](k) \defeq \begin{cases}
		i & \text{ if } k=j\\
		j & \text{ if } k=i\\
		k & \text{ otherwise}
	\end{cases}
\end{align}
It is known that $T$ is generated by the replacements and transpositions, in particular, each non-surjective transformation is the composition of some replacements, and each surjective transformation (i.e., permutation) is a composition of some transpositions.\footnote{See, for example, \cite[Exercise 6(a), p. 39.]{Howie95} and \cite[Remark 1]{Howie78}. Recall that  $\alpha$ is finite.}
 
Let $\c_i$, $\s_{ij}$, and $\p_{ij}$ denote $\c_{(\{ i\})}$, $\s_{[i\slash j]}$, and $\s_{[i,j]}$, respectively. All the other operations are term-definable from these in $\RPA_{\alpha}$. For example, $\c_{(\{ 1,2\})}x=\c_1\c_2x$, $s_{[1\slash2]\circ[0\slash 3]}x=\s_{12}\s_{03}x$ and $\s_{[0,1]\circ[0,3]}x=\p_{01}\p_{03}x$. 
Let $\FRPA_{\alpha}$ denote the class of all reducts of algebras in $\RPA_{\alpha}$ when we keep only $\c_i, \s_{ij}$ and $\p_{ij}$ for $i,j<\alpha$ from the extra Boolean operations and let $\SRd_{cs}\FRPA_{\alpha}$ denote the class of subalgebras of algebras in $\FRPA_{\alpha}$ when we discard the permutation operations $\p_{ij}$. We denote the class of algebras in $\RDF_{\alpha}$ when keeping only the cylindrifications $\c_i, i<\alpha$ from $\c_{(\Gamma)}, \Gamma\subseteq\alpha$ also by $\RDF_{\alpha}$. Analogously to $\PA_{\alpha}$, 
Ildik\'o Sain and Richard J.\ Thompson in \cite[Definition 1]{Sain-Tho} defined $\FPA_{\alpha}$ by a finite set of equations and they proved that $\PA_{\alpha}$ and $\FPA_{\alpha}$ are term-definitionally equivalent.
The equations in the two theorems below are the defining equations of $\FPA_{\alpha}$ in which $\s_{ij}$ and $\p_{ij}$ occur.
\begin{theorem}[$\s_{ij}$ is finitely axiomatizable over $\c_{i}$]\label{main2-thm} Assume that $3\le\alpha$ is finite.
	Let  $\gA=\langle A,+,-,\c_i,\s_{ij}\rangle_{i,j<\alpha}$. Then $\gA\in\SRd_{cs}\FRPA_{\alpha}$ if and only if
	$\langle A,+,-,\c_{i}\rangle_{i<\alpha}\in\RDF_{\alpha}$ and the following equations hold in $\gA$ for all distinct $i,j,k<\alpha$.
	\begin{align}
		&	\s_{ij}(x+y)=\s_{ij}(x)+\s_{ij}(y),\quad \s_{ij}(-x)=-\s_{ij}(x), \label{main6}\\
		&	\s_{ij}\c_ix = \c_ix,\quad \s_{ij}x=\c_i\s_{ij}x,\quad \s_{ii}(x)=x, \label{main7}\\
		&  \s_{ij}\c_kx=\c_k\s_{ij}x,\label{main8} \\
		&	\s_{ij}\s_{kj}x=\s_{ij}\s_{ki}x. \label{main9}			
	\end{align}
\end{theorem}

\begin{theorem}[$\p_{ij}$ is finitely axiomatizable over $\c_{i},\s_{ij}$]\label{main3-thm} Assume that $3\le\alpha$ is finite.
	Let  $\gA=\langle A,+,-,\c_{i},\s_{ij},\p_{ij}\rangle_{i,j<\alpha}$. Then $\gA\in\FRPA_{\alpha}$ if and only if
	$\langle A,+,-,\c_{i},\s_{ij}\rangle_{i,j<\alpha}\in\SRd_{cs}\FRPA_{\alpha}$ and the following equations hold in $\gA$ for all distinct $i,j,k<\alpha$.
	\begin{align}
		&	\p_{ij}(x+y)=\p_{ij}(x)+\p_{ij}(y),\quad \p_{ij}(-x)=-\p_{ij}(x), \label{main10}\\
		&	\p_{ij}\p_{ij}x = x,\quad \p_{ij}x=\p_{ji}x,\quad \p_{ii}x=x, \label{main11}\\
		&	\p_{ij}\s_{ij}x=\s_{ji}\p_{ij}x, \label{main12}\\	
		&  \p_{ij}\p_{ik}x=\p_{jk}\p_{ij}x. \label{main13}
	\end{align}
\end{theorem}
We do not know whether $\p_{ij}$ is finitely axiomatizable over $\c_i$, but we know that $\s_{ij}$ is finitely axiomatisable over $\c_i,\p_{ij}$, see \cite{AGyN-1}. 

The question of which of the extra Boolean operations of $\FRPEA_{\alpha}$ bring in infinitely many new equations is systematically investigated in \cite{AAPALI97ec}. The answer to this question depends on the order in which the operations are added to our algebras; in \cite[Figure 1]{AAPALI97ec} the known answers as of 1991 are indicated on the edges of a tree. 
Theorems~\ref{main2-thm} and \ref{main3-thm} above provide answers for two edges on that tree marked as unknown. 

In Figure \ref{fig:tree}, we reproduce an updated version of this tree where we indicate new information obtained after 1991. Five edges are marked in \cite[Figure 1]{AAPALI97ec} as unknown. Two of them are answered in the present paper, one was answered by Monk \cite{M99} in 1997, and another one in \cite{AGyN-1}. As of now, only one answer remains unknown, namely whether $\p_{ij}$ is finitely axiomatizable over $\c_i$. Some of the unpublished proofs used in \cite[Figure 1]{AAPALI97ec} have appeared in the meantime, these are also indicated in Figure \ref{fig:tree} by labels on the edges. Monk's proof in \cite{M99} is still unpublished, but another proof for the result is given in \cite[Theorem 2.6(ii)]{AGyN-1}.

\begin{figure}[h]\begin{center}
	\begin{tikzpicture}[x=1cm,y=1cm]
	  \node[font=\small] at (1.1,13.2) {$3\le n<\omega$};

	  \node[op] (bool) at (1.15,7.25) {$\cup,\setminus$};

	  \node[op] (dtop) at (3.15,10.25) {$\d_{ij}$};
	  \node[op] (ctop) at (4.55,11.55) {$\c_i$};
	  \node[op] (ptop) at (6.45,12.82) {$\p_{ij}$};

	  \draw[fin] (bool) -- node[edgelab, pos=.43, left=2pt] { } (dtop);
	  \draw[nonfin] (dtop) -- (ctop);
	  \node[smalllab,align=center] at (3.40,11.22) { };
	  \draw[nonfin] (ctop) -- (ptop);
	  \node[smalllab,align=center] at (5.28,12.42) {\cite{ANTnew}};

	  \node[op] (stop1) at (4.66,10.82) {$\s_{ij}$};
	  \draw[fin] (dtop) -- node[edgelab, pos=.58, below, sloped] { } (stop1);

	  \node[op] (ctop2) at (6.55,11.67) {$\c_i$};
	  \draw[nonfin] (stop1) -- node[edgelab, pos=.53, above, sloped] { } (ctop2);

	  \node[op] (ptop2) at (6.52,10.78) {$\p_{ij}$};
	  \draw[fin] (stop1) -- node[edgelab, pos=.52, below] { } (ptop2);

	  \node[op] (ctop3) at (9.12,10.85) {$\c_i$};
	  \draw[nonfin] (ptop2) -- (ctop3);

	  \node[op] (ptop3) at (4.62,9.86) {$\p_{ij}$};
	  \draw[fin] (dtop) -- node[edgelab, pos=.54, below] { } (ptop3);

	  \node[op] (ctop4) at (6.75,10.39) {$\c_i$};
	  \draw[nonfin] (ptop3) -- (ctop4);

	  \node[op] (stop2) at (6.62,9.67) {$\s_{ij}$};
	  \draw[fin] (ptop3) -- node[edgelab, pos=.53, below] { } (stop2);

	  \node[op] (cmid) at (3.15,7.55) {$\c_i$};
	  \draw[nonfin] (bool) -- node[edgelab, pos=.56, above, sloped] { } (cmid);

	  \node[op] (dmid1) at (4.85,8.48) {$\d_{ij}$};
	  \draw[nonfin] (cmid) -- node[edgelab, pos=.55, above, sloped] { } (dmid1);

	  \node[op] (smid) at (4.92,7.55) {$\s_{ij}$};
	  \draw[annfinred] (cmid) -- (smid);

	  \node[op] (pmidL) at (4.90,6.56) {$\p_{ij}$};
	  \draw[unknown] (cmid) -- (pmidL);

	  \node[op] (dmid2) at (6.72,8.45) {$\d_{ij}$};
	  \draw[nonfin] (smid) -- node[edgelab, pos=.53, above, sloped] { } (dmid2);

	  \node[op] (pmidR) at (6.76,7.52) {$\p_{ij}$};
	  \draw[annfinred] (smid) -- (pmidR);

	  \node[op] (dmid3) at (8.42,7.52) {$\d_{ij}$};
	  \draw[nonfin] (pmidR) -- node[edgelab, pos=.53, above] { } (dmid3);

	  \node[op] (dmid4) at (6.82,6.63) {$\d_{ij}$};
	  \draw[nonfin] (pmidL) -- node[edgelab, pos=.52, above, sloped] { } (dmid4);

	  \node[op] (smid2) at (6.86,5.87) {$\s_{ij}$};
	  \draw[fin] (pmidL) -- node[edgelab, pos=.53, below] {\cite{AGyN-1}} (smid2);

	  \node[anntext,text=annred] (rednote) at (10.08,8.05)
	    {Theorem \ref{main2-thm}};
	  \draw[-{Stealth[length=2mm]},draw=annred,line width=.85pt]
	    (rednote.west) .. controls (8.05,8.75) and (6.35,8.95) .. (4.60,7.82);

	  \node[anntext,text=annred] (rednote2) at (10.08,7.05)
	    {Theorem \ref{main3-thm}};
	  \draw[-{Stealth[length=2mm]},draw=annred,line width=.85pt]
	    ($(rednote2.west)+(0,-.10)$) .. controls (8.10,5.85) and (7.28,6.70) .. (6.02,7.48);

	  \node[op] (slow) at (3.15,4.90) {$\s_{ij}$};
	  \draw[fin] (bool) -- node[edgelab, pos=.54, right] {\cite{SN96}} (slow);

	  \node[op] (dlow1) at (4.88,5.26) {$\d_{ij}$};
	  \draw[fin] (slow) -- node[edgelab, pos=.52, above, sloped] { } (dlow1);

	  \node[op] (clow1) at (5.08,4.42) {$\c_i$};
	  \draw[nonfin] (slow) -- (clow1);

	  \node[op] (plow1) at (4.92,3.58) {$\p_{ij}$};
	  \draw[fin] (slow) -- node[edgelab, pos=.53, below, sloped] { } (plow1);

	  \node[op] (dlow2) at (6.75,4.10) {$\d_{ij}$};
	  \draw[fin] (plow1) -- node[edgelab, pos=.52, above, sloped] { } (dlow2);

	  \node[op] (clow2) at (6.93,3.35) {$\c_i$};
	  \draw[nonfin] (plow1) -- (clow2);

	  \node[op] (plow2) at (3.15,1.90) {$\p_{ij}$};
	  \draw[fin] (bool) -- node[edgelab, pos=.55, left=1pt] {} (plow2);

	  \node[op] (dlow3) at (4.86,1.91) {$\d_{ij}$};
	  \draw[fin] (plow2) -- node[edgelab, pos=.53, above] { } (dlow3);

	  \node[op] (clow3) at (5.05,1.15) {$\c_i$};
	  \draw[nonfin] (plow2) -- (clow3);

	  \node[op] (slow2) at (4.94,0.03) {$\s_{ij}$};
	  \draw[fin] (plow2) -- node[edgelab, pos=.55, below, sloped] {} (slow2);

	  \node[smalllab,align=left] (bottomnote) at (8.48,1.35)
	    {\cite{M99} (\cite{AGyN-1})};
	  \draw[-{Stealth[length=2mm]},line width=.85pt]
	    (bottomnote.west) .. controls (6.95,1.62) and (5.95,1.58) .. (4.25,1.53);
	\end{tikzpicture}
\end{center}
\caption{Cartesian space ${}^{n}U$ with operations: those along the path leading to the node;
solid thick lines: not finitely axiomatizable; dashed lines: finitely axiomatizable; solid thin line: unknown (to the authors).}
\label{fig:tree}
\end{figure}

With these two lines filled out, an interesting behaviour of the substitution operations $\s_{ij}$ can be read off the tree: these operations never bring in infinitely many new equations. In  contrast, the permutation operations $\p_{ij}$ exhibit a mixed behaviour: they bring in infinitely many new axioms in the presence of the diagonal constants and the cylindrifications (see \cite{ANTnew}), and when $\s_{ij}$ are available but the diagonals are not, they do not contribute to infinite axiomatizability, by Theorem \ref{main3-thm} above.

The substitution operations $\s_{ij}$ behave rather interestingly in other respects, too, and they are quite well investigated.%
\footnote{See, for example \cite{Fer}, \cite[sections 1.5, 1.11]{HMT}, \cite{P1}, \cite{Pr70}, \cite{SN96}, \cite{Say}, \cite{Thompson1993}.}
For example, Charles Pinter proved in \cite{P1} that in certain complete algebras %
cylindrifications $\c_i$ and diagonals $\d_{ij}$ can be defined from the $\s_{ij}$ so that these defined operations satisfy the defining equations of $\CA_{\alpha}$. This fact will be one of the key steps of our proofs. %
 Because these intimate connections between  $\s_{ij}$ and the $\c_i, \d_{ij}$ hold only in complete algebras, we will rely on the existence of the canonical embedding algebras, because they are complete. Robert Goldblatt's results on canonical varieties come to our aid here. We note that Goldblatt has several papers where he applies the duality between structures and complex algebras to algebras of relations, see \cite{Gol89, Gol91, Goldblatt, AGolN}.

Originally, Halmos defined polyadic algebras for algebraizing first-order logic without equality, while Tarski defined cylindric algebras for algebraizing first-order logic with equality. The equality atomic formulas in first-order logic are algebraized by the diagonal constants $\d_{ij}$. Substituting variables in formulas is treated in polyadic algebras explicitly via the substitution operations $\s_{\sigma}$, while they are treated implicitly via defined terms in cylindric algebras, see \cite[sec.1.5, end of sec.1.11]{HMT}. Halmos \cite[p.28]{Halmos} raises the problem of understanding the relation between these two treatments of equality and variable-substitutions. The results in the present paper contribute to this investigation. Detailed connections between variable-substitution in formulas and the algebraic substitution operations can be found in \cite{ANTnew} and in \cite{AGyN-1}. It turns out that already on the syntactic level of logic, there is a striking difference between transposition of variables ($\p_{ij}$) and replacement of variables ($\s_{ij}$). For more on consequences of these axiomatizability investigations in first-order logic see Monk~\cite{M71}, Johnson~\cite{J73} and \cite[Sections 1,9]{ANTnew}.

The paper is organized as follows. In Section \ref{sec:2}, we recall the definitions of the main classes of algebras relevant to our study. Section \ref{sec:3} shows that our classes are canonical, that is, they are closed under taking canonical embedding algebras. In Section \ref{sec:fpa} we introduce equational axiom systems. Then, in Section \ref{sec:main1}, we prove 
Theorem~\ref{main2-thm}, and in Section~\ref{sec:povercs} we prove Theorems~\ref{main3-thm}, \ref{main1-thm}. Finally, in the Appendix we sketch
alternative proofs for Theorems~\ref{thm:1a} and~\ref{thm:2a}, suggested by the anonymous referee, which avoid atoms and union-completeness.

\section{Setting the stage}\label{sec:2}

In this section, we introduce $\RDF_{\alpha}$ and $\RPA_{\alpha}$, together with some algebras of $\alpha$-ary relations that we will deal with in the present paper.

We will use that an ordinal $\alpha$ is the set of smaller ordinals, so $\alpha=\{ i : i<\alpha\}$, and so a $U$-termed $\alpha$-sequence $s=\langle s_i\rangle_{i<\alpha}$ is a function $s:\alpha\to U$ defined by $s(i)=s_i$ for all $i<\alpha$. 

The set of all $U$-termed $\alpha$-sequences is denoted by ${}^\alpha U$. 
For a sequence $s\in {}^{\alpha}U$, $i <\alpha$, 
and $u\in U$ we denote by $s(i/u)$ the sequence we obtain 
from $s$ by replacing its $i$th value with $u$.  We note that 
$s(i/s(j)) = s\circ[i/j]$, and $s(i/s(j))(j/s(i)) = s\circ [i,j]$ where $\circ$ denotes function composition.
An $\alpha$-ary relation over $U$ is a subset 
$R\subseteq {}^{\alpha}U$. The following unary operations
on $\alpha$-ary relations over $U$ are defined:
\begin{align}
	\C_i^U(R) &\defeq \{ s\in {}^{\alpha}U:\; s(i/u)\in R\text{ for some } u\in U\},\\
	\S_{ij}^U(R) &\defeq \{ s\in {}^{\alpha}U :\; s(i\slash s(j))\in R \},\\
	\P_{ij}^U(R) &= \big\{ s\in {}^{\alpha}U :\; s(i/s(j))(j/s(i))\in R \big\},\\
	\D_{ij}^U &\defeq \{ s\in {}^{\alpha}U:\; s(i)=s(j) \}.
\end{align}
We might omit the superscript $U$ when this is not likely to lead to confusion. 
We note that $\P_{ij}^U(R)$ does not, in fact, depend on $U$ (nevertheless, for clarity, we often indicate it), while $\C_i^U(R)$, $\S_{ij}^U(R)$ and $\D_{ij}^U$  do depend on $U$. The set of all subsets of a set $V$ is denoted by $\Sb(V)$, so 
$\Sb({}^\alpha U)$ is the set of all $\alpha$-ary relations on $U$. 

The $\alpha$-dimensional \emph{full finitary polyadic set algebra over $U$}  is the structure
\begin{align}
	\gA = \< \Sb({}^{\alpha}U), \cup, \setminus, \C_i^U, \S_{ij}^U, \P_{ij}^U\>_{i,j<\alpha}\,.
	\label{setfpa}
\end{align}

For a class $K$ of algebras, $\mathbf{I}K,\mathbf{S}K,\mathbf{P}K,\mathbf{Up}K$ denote the classes of all isomorphic copies, all subalgebras, and all 
isomorphic copies of products and of ultraproducts of elements of $K$, respectively. We will focus on the following classes of algebras
\begin{align}
	\RA_{\alpha}^{csp} &\defeq 
		\mathbf{SI}\big\{ \< \Sb({}^{\alpha}U), \cup, \setminus, \C_i^U, \S_{ij}^U, \P_{ij}^U\>_{i,j<\alpha}:\; U\text{ is a set} \big\}, \\
	\RA_{\alpha}^{cs} &\defeq 
		\mathbf{SI}\big\{ \< \Sb({}^{\alpha}U), \cup, \setminus, \C_i^U, \S_{ij}^U\>_{i,j<\alpha}:\; U\text{ is a set} \big\}, \\
	\RA_{\alpha}^{c} &\defeq 
		\mathbf{SI}\big\{ \< \Sb({}^{\alpha}U), \cup, \setminus, \C_i^U\>_{i<\alpha}:\; U\text{ is a set} \big\}.
\end{align}
The superscripts \emph{csp}, \emph{cs} and \emph{c} refer to the similarity types of the algebras: 
\begin{align*}
	c &\defeq \text{Boolean} + \{\c_i:\; i<\alpha\}, \\
	cs &\defeq \text{Boolean} + \{\c_i, \s_{ij}:\; i,j<\alpha\},  \\
	csp &\defeq \text{Boolean} + \{\c_i, \s_{ij}, \p_{ij}:\; i,j<\alpha\}.	
\end{align*}
If $\gA$ is a $cs$-type algebra, then $\Rd_{c}\gA$ denotes the $c$-type reduct of $\gA$, etc. 

Let $\tau$ denote one of \emph{c}, \emph{cs}, \emph{csp}. It is known that $\RA_{\alpha}^\tau$ is a discriminator class and it is closed under taking ultraproducts, see e.g., \cite[Lemmas  7.1.6 and  7.1.7]{UALbook}. Hence, $\RA_{\alpha}^\tau=\mathbf{SUp}\RA_{\alpha}^\tau$ is a class axiomatizable by a set of universal formulas (see \cite[Thm.2.20]{BurrisSankappanavar1981}) and $\mathbf{SP}\RA_{\alpha}^\tau=\mathbf{HSP}\RA_{\alpha}^\tau$ (see \cite[Thm.2.7.5]{UALbook}) is a class axiomatizable by equations, in other words, $\mathbf{SP}\RA^\tau_{\alpha}$ is the variety generated by $\RA_{\alpha}^\tau$. Subalgebras of reducts of full finitary polyadic set algebras are called \emph{set algebras}, so $\RA_{\alpha}^\tau$ is the  class of the $\tau$-type set algebras, up to isomorphisms%
\footnote{In the following, we sometime omit ``up to isomorphisms" when this is not likely to cause confusion.}. The elements of $\mathbf{SP}\RA_{\alpha}^\tau$ are called \emph{representable algebras}.  These are natural algebras of $\alpha$-ary relations. In the literature, $\RA^{c}_{\alpha}$ is called the class of $\alpha$-dimensional \emph{diagonal-free cylindric set algebras} \cite[5.1]{HMT}, and $\RA^{csp}_{\alpha}$ is the class of $\alpha$-dimensional \emph{finitary polyadic set algebras}.  
The classes $\RDF_{\alpha}$ and $\FRPA_{\alpha}$  of representable diagonal-free and representable finitary polyadic algebras are $\mathbf{SP}\RA_{\alpha}^c$ and  $\mathbf{SP}\RA_{\alpha}^{csp}$, respectively. Recall that for finite dimension $\alpha$, finitary polyadic algebras are term-definitionally equivalent to polyadic algebras as defined by Halmos; let $\RPA_{\alpha}$ denote the  class of the latter algebras.

The goal of this paper is to prove that $\RPA_{\alpha}$ is finitely axiomatizable over $\RDF_{\alpha}$ for finite $\alpha\ge 3$.
We will show this by showing that $\RA_{\alpha}^{cs}$ is finitely axiomatizable over $\RA_{\alpha}^{c}$, and  $\RA_{\alpha}^{csp}$ is finitely axiomatizable over $\RA_{\alpha}^{cs}$, and  then using term-definability as proved in \cite{Sain-Tho}.

\section{Canonicity of $\RA^c_{\alpha}$ and $\RA^{cs}_{\alpha}$}\label{sec:3}

In this section, we use the notion of canonical embedding algebras, 
cf.\ \cite[2.7]{HMT}, \cite[2.7.2, 2.7.3]{HiHo}, \cite[Thms.\ 2.8.9, 2.8.15]{UALbook}. First, we recall the basic definitions 
customized for $c$-type algebras.

\begin{definition}
	A ($c$-type) \emph{atom 
	structure} of dimension $\alpha$ is a first-order 
	relational structure
	$\cS = \<X, T_i\>_{i<\alpha}$, where $T_i\subseteq {}^2X$ for every $i<\alpha$.
\end{definition}\medskip

\begin{definition}
	The \emph{full complex algebra} of an atom structure 
	$\cS = \<X$,  $T_i$ $\>_{i<\alpha}$ is the algebra
	\begin{align}
		\Cm(\cS) \defeq \< \Sb(X), T_i^{*} \>,
	\end{align}
	where for any $A\subseteq X$ we let
	\begin{align}
		T_i^{*}(A) \defeq \{ y\in X:\; \<x,y\>\in T_i\text{ for some } x\in A\}\,.
	\end{align}
\end{definition}
\medskip

\begin{definition}
	Let $\gA = \<A, +, -, \c_i\>_{i<\alpha}$ be a $c$-type algebra, such that
	$\<A, +,-\>$ is a Boolean algebra.  
	The \emph{ultrafilter frame} of $\gA$ is the first-order relational structure
	\begin{align}
		\Uf(\gA) \defeq \< \uf(\gA), T_i\>_{i<\alpha},
	\end{align} 
	where $\uf(\gA)$ is the set of ultrafilters of $\gA$, and $T_i\subseteq {}^{2}\uf(\gA)$
	is defined as
	\begin{align}
		\<p, q\> \in T_i &\quad\Longleftrightarrow\quad \{ \c_i(x):\; x\in p\} \subseteq q.
	\end{align}
\end{definition}\medskip

\begin{definition} $\gA$ is called a \emph{Boolean algebra with operators} ($\BAO$) if the Boolean 
reduct of $\gA$ is a Boolean algebra, and the extra Boolean operations are additive.
 $\gA$ is called a \emph{complete BAO} if $\gA$ is a BAO, arbitrary suprema exist in the Boolean reduct of $\gA$ and the extra Boolean operations preserve infinite sums argumentwise.
  $\gA$ is called a \emph{union-complete BAO} if $\gA$ is a complete BAO with Boolean reduct a set Boolean algebra closed under arbitrary unions. This implies that $\sum X=\bigcup X$ for all nonempty $X$. 
\end{definition} 

The \emph{canonical embedding algebra} of a $\BAO$ $\gA$, 
denoted by $\Em(\gA)$, is the full complex algebra of its ultrafilter frame 
\begin{align}
	\Em(\gA) \defeq \Cm\Uf(\gA) .
\end{align}
By J\'onsson--Tarski \cite{JT}, every $\BAO$
$\gA$ embeds into $\Em(\gA)$, and  $\Em(\gA)$ is a union-complete BAO. Another important property of the canonical embedding algebras $\Em(\gA)$ is that they preserve positive equations. A \emph{positive equation} is an equation in the language of $\gA$ in which negation $-$ does not occur, but the derived operations $\cdot, 0, 1$ can occur. Thus, $\Em(\gA)\models e$ if $\gA\models e$, for any positive equation $e$.  See, e.g., \cite[2.8.9, 2.8.15]{UALbook}, \cite[2.7.5, 2l7.14]{HMT}.

 A class $K$ of algebras is called \emph{canonical} if $\Em(\gA)\in K$ whenever $\gA\in K$.
\bigskip
 
The following theorem is a corollary of Robert Goldblatt's results on the duality between structures and Boolean algebras with operators.

\begin{theorem}\label{thm:canonical}
		$\RA_{\alpha}^\tau$ as well as $\mathbf{SP}\RA_{\alpha}^\tau$ are canonical, for $\tau\in\{ c, cs, csp\}$.
\end{theorem}
\begin{proof}
	We elaborate the proof for $\tau=c$, the proofs for the other cases can be obtained by making obvious modifications. A \emph{c-type sequence-structure} is
	\begin{align}
		\langle{}^\alpha U, T_i\rangle_{i<\alpha}
	\end{align}
	where $U$ is a non-empty set and $sT_iz$ iff $z=s(i\slash u)$ for some $u\in U$. Let $\Str^c$ denote the class of all $c$-type sequence-structures.
	By definition, 
	\begin{align}\label{up}
		\RA^c_{\alpha} = \mathbf{SI}\{ \Cm(\cS) : \cS\in\Str^c\} .
	\end{align}
	Now, $\Str^c$ is closed under taking ultraproducts by \cite[Lemma 3.5]{Goldblatt}. Hence, $\RA^c_{\alpha}$ and $\mathbf{SP}\RA^c_{\alpha}$ are canonical, by Corollary 3.6.3 and Theorem 3.6.7 of \cite{Gol89}.
\end{proof}

All the above, e.g., the definitions of atom structure, full complex algebra, and ultrafilter frame
can be extended in a natural way to $cs$-type  and $csp$-type algebras. For completeness,
we briefly include the definitions for $cs$-type algebras.

\begin{definition}
	A ($cs$-type) atom structure of dimension $\alpha$ is a first-order 
	relational structure
	$\cS = \<X, T_i, R_{ij}\>_{i,j<\alpha}$, where $T_i, R_{ij}\subseteq {}^2X$ 
	for every $i,j<\alpha$.
\end{definition}\medskip

\begin{definition}
	The full complex algebra of a $cs$-type 
	atom structure $\cS = \<X$,  $T_i$, $R_{ij}$ $\>_{i,j<\alpha}$ is the
	algebra
	\begin{align}
		\Cm(\cS) = \< \Sb(X), T_i^{*}, R_{ij}^{*} \>,
	\end{align}
	where for any binary relation $P\subseteq {}^2X$ we define 
	$P^{*}:\Sb(X)\to\Sb(X)$ by
	\begin{align}
		P^{*}(A) = \{ y\in X:\; \<x,y\>\in P\text{ for some } x\in A\}\,.
	\end{align}
\end{definition}
\medskip

\begin{definition}
	Let $\gA = \<A, +, -, \c_i, \s_{ij}\>_{i,j<\alpha}$ be a $cs$-type algebra, 
	such that
	$\<A, +, - \>$ is a Boolean algebra.  
	The \emph{ultrafilter frame} of $\gA$ is the first-order relational structure
	\begin{align}
		\Uf(\gA) = \< \uf(\gA), T_i, R_{ij}\>_{i,j<\alpha},
	\end{align} 
	where $\uf(\gA)$ is the set of ultrafilters of $\gA$, 
	and $T_i, R_{ij}\subseteq {}^{2}\uf(\gA)$
	are defined as
	\begin{align}
		\<p, q\> \in T_i &\quad\Longleftrightarrow\quad 
			\{ \c_i(x):\; x\in p\} \subseteq q, \\
		\<p, q\> \in R_{ij} &\quad\Longleftrightarrow\quad  
			\{ \s_{ij}(x):\; x\in p\}\subseteq q\,.
	\end{align}
\end{definition}\medskip

\section{Finitary polyadic algebras and substitution-cylindric algebras}\label{sec:fpa}
In this section, we introduce the classes $\FPA_{\alpha}$ and $\FSC_{\alpha}$ of algebras defined by finite equational systems.

\begin{definition}[Cf. Definition 1 in \cite{Sain-Tho}]\label{def:FPA}
	By a \emph{finitary polyadic algebra} of dimension $\alpha$,
	$\FPA_{\alpha}$, we mean an algebra 
	$\gA = \<A, +, -, \c_i, \s_{ij}, \p_{ij}\>_{i,j<\alpha}$ in
	which equations (F0-F9) below are valid for every $i, j, k < \alpha$.
	\begin{enumerate}[label=(F\arabic*),start=0]
		\item $\<A, +, -\>$ is a Boolean algebra\footnote{The derived Boolean operations $\cdot$, $0$ etc. are taken as defined in the usual way.}, $\s_{ii}(x)=\p_{ii}(x)=x$, $\p_{ij}(x)=\p_{ji}(x)$. \label{FPA:0}
		\item $x\leq \c_ix$ \label{A:c-extensive}
		\item $\c_i(x+y)=\c_ix+\c_iy$ \label{A:ci-additive}
		\item $\s_{ij}\c_ix = \c_ix$ \label{FPA:3}
		\item $\c_i\s_{ij}x = \s_{ij}x$ if $i\neq j$
		\item $\s_{ij}\c_kx = \c_k\s_{ij}x$ if $k\notin \{i,j\}$
		\item $\s_{ij}$ and $\p_{ij}$ are Boolean endomorphisms 
				\label{FPA:endo}\label{A:s-endomorphism}\label{A:S1}
		\item $\p_{ij}\p_{ij}x=x$ \label{FPA:pij-invertible}
		\item $\p_{ij}\p_{ik}x=\p_{jk}\p_{ij}x$ if $i,j,k$ are all distinct
		\item $\p_{ij}\s_{ji}x=\s_{ij}x$ \label{FPA:9}
	\end{enumerate}
\end{definition}

\noindent 
It is routine to check that full finitary polyadic set algebras \eqref{setfpa} satisfy \ref{FPA:0}-\ref{FPA:9}, so $\RA_{\alpha}^{csp}\subseteq\FPA_{\alpha}$. 

\begin{remark}\label{pozitivak}
	The identities \ref{FPA:0}-\ref{FPA:9} are positive equations. 
	It is only \ref{FPA:endo} which needs a comment.
	The equation $-\s_{ij}(x) = \s_{ij}(-x)$ 
	is part of what	it means that $\s_{ij}$ is a Boolean 
	endomorphism, and it is visible that the negation 
	symbol is used here. This equation, as noted in \cite[p.561]{Sain-Tho}, 
	can be replaced by the ones 
	expressing that $\s_{ij}$ preserves $+$, $\cdot$, $0$ and $1$. 
	This is because if $\s_{ij}$ preserves $+$, $\cdot$, $0$ and $1$, then
	\begin{align}
		\s_{ij}(x)\cdot\s_{ij}(-x) &= \s_{ij}(x\cdot -x) = \s_{ij}(0) = 0,\\
		\s_{ij}(x)+\s_{ij}(-x) &= \s_{ij}(x + -x) = \s_{ij}(1) = 1,			
	\end{align}
	and these imply $\s_{ij}(-x) = -\s_{ij}(x)$. The same comment applies to
	the case of $\p_{ij}$. Hence, $\FPA_{\alpha}$ is a canonical variety.
\end{remark}

In Section \ref{sec:povercs}, we will use consequences \ref{uj:p2}--\ref{uj:p4-2} of the $\FPA_{\alpha}$-axioms. In Definition \ref{def:ax-alpha} below and in section~\ref{sec:main1} we will use consequence (FS) of the $\FPA_{\alpha}$-axioms that concerns solely the substitution operations $\s_{ij}, i,j<\alpha$. It is proved in \cite[Remark 7]{Sain-Tho} that (FS) does not follow from the axioms of $\FPA_{\alpha}$ that do not contain any reference to $\p_{ij}, i,j<\alpha$.
	
\begin{lemma}\label{pc-lemma}\label{lem:ST-1.6}
	The following identities hold in $\FPA_{\alpha}$.
	\begin{enumerate}[label=(P\arabic*),start=1]
		\item $\p_{ij}\c_ix = \s_{ji}\c_ix$ \label{uj:p2} 
		\item $\p_{ij}\c_jx = \s_{ij}\c_jx$ \label{uj:p3} 
		\item $\p_{ij}\c_kx = \s_{kj}\s_{ji}\s_{ik}\c_kx$ if $k\notin\{ i,j \}$
		 \label{uj:p1} 
		 \item $\p_{ij}\c_ix = \c_j\p_{ij}x$ for any $i,j<\alpha$ \label{uj:p4-1}
		 \item $\p_{ij}\c_mx = \c_m\p_{ij}x$ if $i,j,m<\alpha$ are distinct. \label{uj:p4-2}
		 \item[(FS)] $\s_{ij}\s_{kj}x=\s_{ij}\s_{ki}x$ if $i,j,k<\alpha$ are distinct. \label{FS}
	\end{enumerate}
\end{lemma}
\begin{proof}
	Here are the derivations.\\

\noindent \ref{uj:p2} \quad $\p_{ij}\c_i \overset{\ref{FPA:3}}{=} \p_{ij}\s_{ij}\c_i 
			\overset{\ref{FPA:0}}{=} \p_{ji}\s_{ij}\c_i 
			\overset{\ref{FPA:9}}{=} \s_{ji}\c_i$
			
\noindent \ref{uj:p3} \quad $\p_{ij}\c_j \overset{\ref{FPA:3}}{=} \p_{ij}\s_{ji}\c_j 
			\overset{\ref{FPA:9}}{=} \s_{ij}\c_j$
			
\noindent \ref{uj:p1} \quad $\p_{ij}\c_k \overset{\ref{FPA:3}}{=} \p_{ij}\s_{ki}\c_k
			\overset{\text{\cite[(15) p.553]{Sain-Tho}}}{=} \s_{kj}\s_{ji}\s_{ik}\c_k$\\
	
\noindent	
\ref{uj:p4-1} and \ref{uj:p4-2}
	are proved as Lemma 1.6 in \cite{Sain-Tho}.%
	\footnote{We note that there are two typos in the statement of Lemma 1.6 in \cite{Sain-Tho}. Namely, 
	(i) writes $\p_{ij}\c_i = \c_i\p_{ij}$, but the correct version is our \ref{uj:p4-1}. In (iii),
	the correct assumption is that $\sigma(i)=i$ and for $j\neq i$, $\sigma(j)\neq i$.
	However, the proof of (the typo-free version of) 
	Lemma 1.6 in \cite{Sain-Tho} is correct, and throughout 
	that paper the lemma is used with the correct statement.}
	(FS) is proved as Lemma 1.5(iii) in \cite{Sain-Tho}.
\end{proof}

Recall that permutations $\sigma:\alpha\to\alpha$ are generated by 
the transpositions. We let $\p_{[i,j]} \defeq \p_{ij}$ and 
$\p_{\sigma\circ \tau} \defeq \p_{\sigma}\circ\p_{\tau}$. This way to 
every permutation $\sigma$ we assign a 
corresponding polyadic term $\p_{\sigma}$. \cite[Claim 1.2]{Sain-Tho} 
shows that different decompositions of the same permutation into 
products of transpositions yield identical terms (in $\FPA_{\alpha}$). 
For a precise definition we refer to \cite[p.554]{Sain-Tho}.
This leads to the important lemma below that we make use of
many times in Section \ref{sec:povercs}.

\begin{lemma}\label{lem:permulem}
	The following identity is a consequence of the $\FPA_{\alpha}$ axioms
	for any permutation $\tau:\alpha\to\alpha$ and $i<\alpha$. 
	\begin{enumerate}[label=(P\arabic*),start=6]
		\item $\p_{\tau}\c_ix = \c_{\tau(i)}\p_{\tau}x$\label{lemma1.6-kov}
	\end{enumerate}
\end{lemma}
\begin{proof}
	Apply Lemma \ref{lem:ST-1.6} iteratively. 
\end{proof}


Next we turn to algebras without the $\p_{ij}$ operations. The
class $\FSC_{\alpha}$
of algebras is axiomatized by the 
finitary polyadic algebra axioms of $\FPA_{\alpha}$ 
\emph{without the $\p_{ij}$'s}, that is, when the 
$\FPA_{\alpha}$ axioms speaking about $\p_{ij}$ are deleted; with (FS) added as axiom.  We will see that these algebras are obtained from cylindric algebras by replacing the diagonal constants $\d_{ij}$ with the derived operations $\s_{ij}$. The elements of $\FSC_{\alpha}$ will be called	\emph{substitution-cylindric} algebras. This class was introduced by Pinter~\cite{P2} where it is denoted by $\QA_{\alpha}$, and its members are called quantifier-algebras.

\begin{definition}[Cf.\ Definition~2.1 in \cite{P2}]\label{def:ax-alpha}
	The class $\FSC_{\alpha}$ consists of the 
	algebras $\gA = \< A, +, -, \c_i, \s_{ij}\>_{i,j<\alpha}$ in which
	the equations below are valid for every $i,j,k<\alpha$.
	\begin{enumerate}[label=(F\arabic*),start=0]
		\item $\<A, +, -\>$ is a Boolean algebra, $\s_{ii}(x)=x$ \label{F:0}
		\item $x\leq \c_ix$
		\item $\c_i(x+y)=\c_ix+\c_iy$
		\item $\s_{ij}\c_ix = \c_ix$ \label{A:sijci=ci}
		\item $\c_i\s_{ij}x = \s_{ij}x$ if $i\neq j$ \label{A:cisij=sij}
		\item $\s_{ij}\c_kx = \c_k\s_{ij}x$ if $k\notin \{i,j\}$ 
		\label{F:beforelast}\label{A:sijck=cksij}
		\item $\s_{ij}$ is a Boolean endomorphisms \label{F:last}
		 \item[(FS)] $\s_{ij}\s_{kj}x=\s_{ij}\s_{ki}x$ if $i,j,k<\alpha$ are distinct.
	\end{enumerate}
	We deliberately use the same labels as in Definition \ref{def:FPA}, even 
	if (F0) and (F6) here are not exactly the same as (F0) and (F6) in Definition
	\ref{def:FPA} as here the $\p_{ij}$'s are missing. We believe that this
	slight abuse of notation will not cause confusion.
\end{definition}

\noindent Clearly, Remark \ref{pozitivak} applies to the $\FSC_{\alpha}$ axioms, 
yielding that these identities are positive, and $\FSC_{\alpha}$ is a canonical variety.
Also, $\RA_{\alpha}^{cs}\subseteq\FSC_{\alpha}$. \\

In the upcoming proofs we will extensively use that diagonal constants $\d_{ij}$ can be adjoined to complete $\FSC_{\alpha}$ algebras so that the cylindric axioms hold. First we recall the definition of a cylindric algebra. 

\begin{definition}[Definition 1.1.1 in \cite{HMT}]\label{def:cax-alpha}
	The class $\CA_{\alpha}$ of \emph{cylindric algebras} consists of the 
	algebras $\gA = \< A, +, -, \c_i, \d_{ij}\>_{i,j<\alpha}$ in which
	the equations below are valid for every $i,j,k<\alpha$.
	\begin{enumerate}[label=(C\arabic*),start=0]
		\item $\<A, +, -\>$ is a Boolean algebra
		\item $\c_i0=0$
		\item $x\le \c_ix$
		\item $\c_i(x\cdot\c_iy)=\c_ix\cdot\c_iy$ 
		\item $\c_i\c_jx=\c_j\c_ix$  
		\item $\d_{ii}=1$ 
		\item $\d_{ij}=\c_k(\d_{ik}\cdot\d_{kj})$ if $k\ne i,j$		
		\item $\c_i(\d_{ij}\cdot x)\cdot \c_i(\d_{ij}\cdot-x)=0$ if $i\ne j$ \label{C:last}
	\end{enumerate}
\end{definition}

Let $\gA=\langle A,+,-,\c_i,\d_{ij}\rangle_{i,j<\alpha}\in\CA_{\alpha}$. The \emph{substitution-cylindric reduct} $\Rd_{sc}\gA$ of $\gA$ is $\langle A,+,-,\c_i,\s_{ij}^{\gA}\rangle_{i,j<\alpha}$ where $\s_{ij}^\gA :A\to A$ is the function defined by $\s_{ij}^\gA(x)=\c_i(\d_{ij}\cdot x)$ when $i\ne j$ and $\s_{ii}^\gA(x)=x$ for $x\in A$, and $\Rd_{sc}\CA_{\alpha}=\{ \Rd_{sc}\gA : \gA\in\CA_{\alpha}\}$. The following theorem says that the variety $\FSC_{\alpha}$ is closely related to the substitution-cylindrification reducts of cylindric algebras.

\begin{theorem}\label{rdsc-thm} Let $3\le\alpha$ be finite. Then (i)-(ii) below hold.
	\begin{enumerate}[label=(\roman*)]
		\item $\gA\in\Rd_{sc}\CA_{\alpha}$ for every complete $\gA\in\FSC_{\alpha}$. \label{rdsc1}
		\item $\FSC_{\alpha} = \mathbf{S}\Rd_{sc}\CA_{\alpha}$. \label{rdsc2}
	
	\end{enumerate}
\end{theorem}
\begin{proof}
The proof of \ref{rdsc1} relies on results in Pinter \cite{P1} and Sain-Thompson \cite{Sain-Tho}. 
First we prove that (F0)-(F6),(FS) imply substitution axioms (S1)-(S6) in \cite{P1}.
Indeed, Proposition 6 in \cite{Sain-Tho} proves that the equations 
	\begin{enumerate}[label=(A\arabic*),start=1]
		\item $\s_{ij}\s_{ik}x = \s_{ik}x$, if $i\neq k$ \label{A:last}
		\label{Pinter:S5}
		\item $\s_{ij}\s_{kl}x = \s_{kl}\s_{ij}x$ if $i\notin \{k,l\}$, 
		$k\neq j$  \label{A:S6}\label{Pinter:S6}
	\end{enumerate}
hold in $\FSC_{\alpha}$ and (S1)-(S6) are \ref{F:last}, \ref{F:0}, (FS),
\ref{Pinter:S5}, and \ref{Pinter:S6}. Also, \cite{Sain-Tho} proves 
\[ \mbox{$\c_ix$ is the least element of $\{ y\in\ran\s_{ij}: y\geq x\}$}   \] 
for $i\ne j$ as (4) on p.545, this is condition ($\pi1$) + (B) in \cite{P1}. Let us define $\d_{ij}^*$, for $i\ne j$, as
the infimum of $\{ y : \s_{ij}y=1\}$ in $\gA$, this infimum exists because $\gA$ is complete. Then (32) on p.561 in \cite{Sain-Tho} proves that $\s_{ij}\d_{ij}^*=1$, so
\[ \mbox{$\d_{ij}^*$ is the least element of $\{  y : \s_{ij}y=1\}$}   \] 
which is condition ($\pi2$) + (C) in \cite{P1}. Finally, condition ($\pi3$) in \cite{P1} is exactly our (F5). Thus the hypothesis part  (S1)-(S6), ($\pi1$)-($\pi3$) and (B),(C) of \cite[Theorem 2.7(ii)]{P1} hold.%
\footnote{In fact, $\FSC_{\alpha}$ can be axiomatized by (S1)-(S6),(B),$(\pi3)$.}
 Then $\gA^*=\langle A,+,-,\c_i,\d_{ij}^*\rangle_{i,j<\alpha}\in\CA_{\alpha}$ in which $\s_{ij}(x)=\c_i(\d_{ij}^*\cdot x)$ holds for $i\ne j$, by \cite[Theorem 2.7(ii)]{P1}. Now, $\gA=\Rd_{sc}\gA^*$ and $\gA^*\in\CA_{\alpha}$ witness that $\gA\in\Rd_{sc}\CA_{\alpha}$, and this finishes the proof of (i). 

To prove (ii), notice first that $\FSC_{\alpha}\subseteq\mathbf{S}\Rd_{sc}\CA_{\alpha}$ follows from (i) together with $\FSC_{\alpha}$ being a canonical variety. For the other direction, we have to prove that (F0)-(F6),(FS) hold in $\Rd_{sc}\CA_{\alpha}$. Now, (F0) is (C0)+\cite[1.5.1]{HMT}, (F1) is (C2), and (F2)-(F6),(FS) are \cite[1.2.6(ii), 1.5.8(i), 1.5.9(ii), 1.5.8(ii), 1.5.3, 1.5.10(ii)]{HMT}, respectively.
\end{proof}

\section{$\RA_{\alpha}^{cs}$ is finitely axiomatizable over $\RA^{c}_{\alpha}$}\label{sec:main1}

In this section, we prove our first main result that states that a substitution-cylindric algebra of finite dimension $\ge 3$ is representable iff its substitution-free reduct is representable. A corollary to this 
result is that substitutions over the Booleans together with cylindrifications are finitely axiomatizable (by the axioms of substitution-cylindric algebras).

\begin{theorem}\label{thm:1a}
		For finite $\alpha>2$, an $\alpha$-dimensional algebra 
		$\gA\in\FSC_{\alpha}$ is 
		representable if and only if 
		its sub\-sti\-tu\-ti\-on-free reduct is representable.
\end{theorem}

\begin{proof}
	The left-to-right direction is immediate. For the converse, 
	we show first that it is enough to prove that if the substitution-free reduct of a
	union-complete $\FSC_{\alpha}$ is a $c$-type set algebra, then the algebra itself 
	is isomorphic to a $cs$-type set algebra.

	\begin{lemma}\label{lem:1}
		To prove Theorem \ref{thm:1a} it is enough to show that the statement
		\begin{align*}
			&\text{if }\gA\in\FSC_{\alpha} \text{ is union-complete, and }  
			 \Rd_{c}\gA = \<A,\cup,\setminus,\C_i^U\>_{i<\alpha}, \\ 
			& \text{ for some } U \text{ and } 
			A\subseteq\Sb({}^{\alpha}U), \text{ then } \gA\in\RA^{cs}_{\alpha}
			\tag{$\sharp$}
			\label{*}
		\end{align*}	
		holds.
	\end{lemma}
	\begin{proof}[of Lemma \ref{lem:1}]
		Assume that \eqref{*} holds. 
		Take an arbitrary $\gA  = 
		\<A$, $ +, -$, $\c_i$, $ \s_{ij}$ $\>_{i,j<\alpha}$ 
		$\in \FSC_{\alpha}$, and assume
		that $\Rd_{c}\gA = \<A, +, -, \c_i\>_{i<\alpha}$ is representable, say 
		\begin{align}
			f: \<A, +, -, \c_i\>_{i<\alpha} \quad\rightarrowtail\quad
			\Pi_{k\in K}\<\Sb({}^{\alpha}U_k), \cup, \setminus, \C_i\>_{i<\alpha}
		\end{align}
		holds for some $f$ and $U_k$, $k\in K$. We have to show that $\gA$ is representable, as well.

		Let $k\in K$, denote by $\pi_k$ the projection function, and let 
		\begin{align}
			g:\Rd_{c}\gA\rightarrow 
			\<\Sb({}^{\alpha}U_k), \cup, \setminus, \C_i\>_{i<\alpha}
		\end{align}
		be the composition of $f$ with $\pi_k$. Let $\gB_k$ be the $g$-image of
		$\gA$ and let $i,j<\alpha$. 
		We show that \ref{F:0}-\ref{F:last} (in Definition
		\ref{def:ax-alpha} of $\FSC_{\alpha}$) imply that 
		the $g$-image of the function $\s_{ij}$ on $A$ is a function on $B_k$.  
		To this end, we have to show that if $g(a)=g(b)$ then
		$g(\s_{ij}a)=g(\s_{ij}b)$, for any $a,b\in A$. Note that $g$ is a
		Boolean homomorphism that preserves also the cylindrifications $\c_i$.
		Now, $g(a)=g(b)$ iff, by $g$ being a Boolean homomorphism,
		$g(a\oplus b)=0$, where $\oplus$ denotes the Boolean symmetric
		difference. Similarly, $g(\s_{ij}a)=g(\s_{ij}b)$ iff
		$g(\s_{ij}a\oplus\s_{ij}b)=g(\s_{ij}(a\oplus b))=0$,
		by \ref{F:last}. Therefore, it is enough to show that
		$g(\s_{ij}x)=0$ whenever $g(x)=0$. Indeed,
		\begin{align*}
			g(x) = 0 
				&\quad\Longrightarrow\quad g(\c_ix)=\c_ig(x)= 0 
				\tag{by $g$ a $\c_i$-homomorphism + $\C_i\emptyset=\emptyset$} \\ 
				&\quad\Longrightarrow\quad g(\s_{ij}\c_ix)=0 \tag{by \ref{A:sijci=ci}} \\
				&\quad\Longrightarrow\quad g(\s_{ij}x)\le g(\s_{ij}\c_ix)=0 
				\tag{by \ref{A:c-extensive}+\ref{A:s-endomorphism}+ $g$ homomorphism}   \\
				&\quad\Longrightarrow\quad g(\s_{ij}x)=0  
		\end{align*}
		\noindent
		We have seen that the $g$-images of the (abstract) substitution operations of $\gA$ 
		are functions on $\gB_k$. 
		Let $\s_{ij}^{+}$ denote these projected functions and let 
		$\gB_k^+ = \< \gB_k, \s_{ij}^{+}\>_{i,j<\alpha}$, for all $k\in K$. 
		Then $g:\gA\to\gB_k^+$ and so $\gB_k^+\in\FSC_\alpha$ by $\gA\in\FSC_\alpha$. 
		Also $\Rd_{c}\gB_k^+ = \gB_k\in\RA_{\alpha}^{c}$ by its construction.

		Consider now the canonical embedding algebra $\Em(\gB_k^+)$. By
		J\'onsson--Tarski \cite{JT} (cf.\ \cite[Thm.~2.8.9]{UALbook},
		\cite[2.7.13]{HMT}), $\gB_k^+$ embeds into $\Em(\gB_k^+)$ and
		$\Em(\gB_k^+)$ is a union-complete BAO. Positive equations are
		preserved by canonical extensions, so Remark~\ref{pozitivak} gives
		$\Em(\gB_k^+)\in\FSC_{\alpha}$.
		Moreover,
		\begin{align*}
			\Rd_c\Em(\gB_k^+)
			=\Em(\Rd_c\gB_k^+)
			=\Em(\gB_k).
		\end{align*}
		By Theorem \ref{thm:canonical}, as $\gB_k\in\RA^c_{\alpha}$, we have
		$\Em(\gB_k)\in\RA^c_{\alpha}$. That is, $\Em(\gB_k)$ is isomorphic to
		a $c$-type set algebra. This isomorphism is a complete isomorphism: For cylindric algebras 
		this is 
		\cite[Theorem 3.4.3]{HiHoComp}, the $c$-type version follows from this and
		\cite[Theorem 3.36]{HiHo}, by adjoining the set-theoretic diagonal 
		elements to the set representation. Identify $\Em(\gB_k^+)$ with this
		union-complete isomorphic copy.
		Then \eqref{*} applies to
		$\Em(\gB_k^+)$, and hence $\Em(\gB_k^+)\in\RA^{cs}_{\alpha}$.
		Since $\gB_k^+$ is a subalgebra of $\Em(\gB_k^+)$, it follows that
		$\gB_k^+\in\RA^{cs}_{\alpha}$.
		This holds for every $k\in K$ and so $\gA\in\mathbf{SP}\RA_{\alpha}^{cs}$ 
		since $f:\gA\rightarrowtail \Pi_{k\in K}\gB_k^+$.
		Consequently, $\gA$ is representable.
	\end{proof}	
		
	We have seen that it is enough to prove \eqref{*}. To do so, take a union-complete
	$\gA\in\FSC_{\alpha}$ and base set $U$ and $A\subseteq\Sb({}^{\alpha}U)$,
	such that
	\begin{align}
		\Rd_c\gA=\<A,\cup,\setminus,\C_i^U\>_{i<\alpha}.
	\end{align}
	Here $A$ is closed under arbitrary
	unions and $\sum X=\bigcup X$ for every $X\subseteq A$, and the additional operations 
	are completely additive.
	Putting back the substitution operations, we write
	\begin{align}
		\gA=\<A,\cup,\setminus,\C_i^U,\s_{ij}^{\gA}\>_{i,j<\alpha},
	\end{align}
	where the operations $\s_{ij}^{\gA}$ are abstract, not yet
	set-represented operations.

By Theorem \ref{rdsc-thm}\ref{rdsc1} then $\gA=\Rd_{sc}\gB$ for some $\gB\in\CA_{\alpha}$, since $\gA\in\FSC_{\alpha}$ is complete. 
By the definition of $\Rd_{sc}\gB$ this means that there are $\d_{ij}^\gA\in A$ for $i,j<\alpha$ such that 
	$\<A, \cup, \setminus, \C_i^U, \d_{ij}^{\gA}\>_{i,j<\alpha}$ is a 
cylindric algebra, where
\begin{align}
	\s_{ij}^{\gA}x = \begin{cases}
		x & \mbox{ if } i=j,\\
		\C_i^U(x\cap \d_{ij}^{\gA}) & \mbox{ if } i\neq j
	\end{cases} \label{eq:s=cd}
\end{align}
holds. \\

	 We now change the  representation of $\gA$ in order to arrange that these ``abstract" diagonal constants $\d_{ij}^{\gA}$ be represented close to the set diagonal constants $\D_{ij}^U$. 
	 We cannot achieve $\d_{ij}=\D_{ij}$ because the diagonals $\D_{ij}$ are not finitely axiomatizable over $\RA^{cs}$, see \cite[Theorem 6]{AAPALI97ec}. 
	 We will be able to achieve $\D_{ij}\subseteq \d_{ij}$ and this will be  enough for representing the $\s_{ij}$ well. The idea for achieving $\D_{ij}\subseteq\d_{ij}$ is that first we ``blow up" each element of $U$ and then we will select a bijection below each $\d_{ij}$ along which we will ``rearrange" the sequences so that we get $\D_{ij}\subseteq \d_{ij}$.
	
	We now change the representation by blowing up each element of $U$. Let $W$ be an infinite set such that $|W|>|U|$. The new base set of 
	the representation will be $V\defeq U\times W$. For $s\in{}^{\alpha}U$ and $w\in{}^{\alpha}W$ let $s\times w\in{}^{\alpha}V$ be defined as
	\[ s\times w\defeq\langle (s_i,w_i) : i<\alpha\rangle.\]
	For each sequence $s\in {}^{\alpha}U$ let us write 
	\begin{align}
		\hat{s} \defeq \{ s\times w :\;  w\in{}^\alpha W \}, \label{eq:hats}
	\end{align}
	and for $X\in\Sb({}^{\alpha}U)$ let
	\begin{align}
		F(X) \defeq \bigcup\big\{ \hat{s}:\; s\in X \big\}\in\Sb({}^{\alpha}V) .
	\end{align}
	It is straightforward to check that $F$ is an embedding 
	\begin{align}
		\< A, \cup, \setminus, \C_i^{U} \>_{i<\alpha}
		 \quad\rightarrowtail\quad
		\< \Sb({}^{\alpha}V), \cup, \setminus, \C_i^{V} \>_{i<\alpha}
	\end{align}
	that takes infinite unions to unions, i.e., $F(\bigcup X)=\bigcup\{ F(x) : x\in X\}$ for all $X\subseteq A$.
	Consequently, there is a subset $A'\subseteq \Sb({}^{\alpha}V)$ 
	such that $\<\gA, \d_{ij}^{\gA}\>$ is isomorphic to 
	an algebra
	\begin{align}
		\gA' = 
		\< A', \cup, \setminus, \C_i^{V}, 
		\s_{ij}', \d_{ij}' \>_{i<\alpha},
	\end{align}
	where $\s_{ij}'$ and $\d_{ij}'$ are the isomorphic 
	images of $\s_{ij}^{\gA}$ and $\d_{ij}^{\gA}$, respectively. Therefore, $\gA'$ is a union-complete BAO, with the corresponding reducts being a cylindric algebra and a substitution-cylindric algebra, with $\s_{ij}'$ defined as in \eqref{eq:s=cd}. \\
	
	We now make preparations for changing the representation once more in order to achieve $\D_{ij}\subseteq\d_{ij}$.
		Elements of $\Sb({}^{\alpha}V)$ are $\alpha$-ary relations. For 
	such a relation
	$R\subseteq {}^{\alpha}V$ and $i<\alpha$ we define
	\begin{align}
		\dom_i(R) \defeq \{ s(i):\; s\in R\}\ \subseteq\ V.
	\end{align} 
	In what follows, $At$ is the set of atoms of $\gA'$.
	\begin{lemma}\label{dom-particio}
		Let $i\neq j<\alpha$, and let $a, b\subseteq \d_{ij}'$ be two atoms. 
		Then either
		\begin{align}
			\dom_i(a) = \dom_i(b),\quad\text{ and }\quad \dom_j(a) = \dom_j(b),
			\label{dom-particio-1}
		\end{align}
		or
		\begin{align}
			\dom_i(a) \cap \dom_i(b) = \emptyset,\quad\text{ and } \quad
			\dom_j(a) \cap \dom_j(b) = \emptyset.\label{dom-particio-2}
		\end{align}
	\end{lemma}
	\begin{proof}[of Lemma \ref{dom-particio}]
		Let $a, b$ be distinct atoms below $\d_{ij}'$, and 
		set $\Gamma=\alpha-\{i,j\}$. Below, $\C_{\Gamma}^V(x)$ denotes $\C_{i1}^V\dots\C_{ig}^V(x)$ where $\Gamma=\{ i1,\dots,ig\}$. This notation makes sense 
		by the properties of $\C_i$.
		By \cite[1.3.3]{HMT}, 
		$\C_{\Gamma}^V(\d_{ij}') = \d_{ij}'$, and
		by \cite[1.10.3]{HMT}, $\C_{\Gamma}^V(a)$ is an atom
		of the algebra of $\Gamma$-closed elements. 
		In particular, $\C_{\Gamma}^V(a)$ and $\C_{\Gamma}^V(b)$ are 
		below $\d_{ij}'$, and are either equal or disjoint.
		By definition of $\C_{\Gamma}^V$, we have
		\begin{align}
			\dom_i(a) = \dom_i(\C_{\Gamma}^V(a)), \quad
			\dom_j(a) = \dom_j(\C_{\Gamma}^V(a)),  \\
			\dom_i(b) = \dom_i(\C_{\Gamma}^V(b)), \quad
			\dom_j(b) = \dom_j(\C_{\Gamma}^V(b))\,.
		\end{align}
		If $\C_{\Gamma}^V(a) = \C_{\Gamma}^V(b)$, then 
		\eqref{dom-particio-1} follows. In the other case, 
		assume that 
		\begin{align}
			\C_{\Gamma}^V(a) \cap \C_{\Gamma}^V(b) = \emptyset\,. \label{secondcase}
		\end{align}
		We show \eqref{dom-particio-2}.
		By way of contradiction, suppose there is $v\in \dom_{i}(a)\cap \dom_i(b)$.
		Then there are $x,y\in {}^{\alpha}V$ such that $x\in a$, 
		$y\in b$, and $x(i)=y(i)=v$. Let $z = x(j/y(j))$. 
		Then $z\in \C_{\Gamma}^V(b)$, because $z$ and $y$ can differ only
		at their $\Gamma$-coordinates. Also, $z\in \C_j^V(a)$, because $x$ and $z$
		may differ only at the $j$th coordinate. But then
		\begin{align}
			z &\in \d_{ij}'\cap -\C_{\Gamma}^V(a)
				\tag{by \eqref{secondcase} and $z\in \C_{\Gamma}^V(b)\subseteq \d_{ij}'$}\\ 
			z &\in \C_j^V(\d_{ij}'\cap \C_{\Gamma}^V(a)) \tag{by 
			$\C_{\Gamma}^V(a)\subseteq \d_{ij}'$, monotonicity of 
			$\C_{\Gamma}^V$ and $z\in \C_j^V(a)$}
		\end{align}
		Using monotonicity of $\C_j^V$ we end up at
		\begin{align}
			z\in \C_j^V(\d_{ij}'\cap \C_{\Gamma}^V(a))\cap 
			\C_j^V(\d_{ij}'\cap -\C_{\Gamma}^V(a)),
		\end{align}
		which contradicts cylindric axiom (C7) (see Definition \ref{def:cax-alpha}):
		\begin{align}
			\CA_{\alpha}\models \c_j(\d_{ij}\cdot x)\cdot\c_j(\d_{ij}\cdot -x) = 0.
		\end{align}
		We proved that if \eqref{secondcase} holds, then 
		$\dom_{i}(a)\cap \dom_i(b)=\emptyset$. Exchanging $i$ and $j$ in the
		argument yields $\dom_{j}(a)\cap \dom_j(b)=\emptyset$.
	\end{proof}
	
	\begin{lemma}\label{egeszV}
		For all $i,j<\alpha$, $\dom_{i}(\d_{ij}') = \dom_{j}(\d_{ij}') = V$.
	\end{lemma}
	\begin{proof}[of Lemma \ref{egeszV}]
		By the cylindric axioms	(see \cite[1.3.2]{HMT}),
		\begin{align}
			\C_j^V(\d_{ij}') = {}^{\alpha}V\,. \label{eq:ujcd1}
		\end{align}
		Now, if for some $v\in V$, $v\notin \dom_{i}(\d_{ij}')$, then
		$v\notin \dom_{i}\big(\C_j^V(\d_{ij}')\big)$, as the $i$th
		coordinates of the sequences in $\d_{ij}'$ and $\C_j^V(\d_{ij}')$
		are the same. But then \eqref{eq:ujcd1} would be impossible. Similarly for $j$.
	\end{proof}
		
	Fix $i\neq 0$
	and let $a\subseteq \d_{0i}'$ be an atom. 
	By construction (of the embedding $F$), 
	if $\<u,w\>\in \dom_{0}(a)$, then $\<u,v\>\in \dom_{0}(a)$ for 
	every $v\in W$, and similarly with $\dom_{i}(a)$. This means
	that every $\<u,w\>\in \dom_{0}(a)$ has $|W|$-many $a$-neighbours in 
	$\dom_{i}(a)$, and vice-versa.
	Using transfinite recursion, we construct a bijective function
	\begin{align}
		f_{0i}^a \subseteq \dom_{0i}(a) \overset{\text{def}}{=} \{ \< s(0), s(i)\>: s\in a\}.
	\end{align}
	To do so, enumerate $\dom_{0}(a)$ and $\dom_{i}(a)$ in $|W|$-type, 
	and starting with 
	the empty list of pairs, 
	in consecutive rounds for each not yet selected element in $\dom_{0}(a)$
	pick one 
	of its $a$-neighbours in $\dom_{i}(a)$ that has not yet 
	been selected; and for each not yet selected element in $\dom_{i}(a)$ 
	find a pair 
	among its $a$-neighbours in $\dom_{0}(a)$ that has not yet been picked. 
	This is possible, because each element has $|W|$ possible pairs in $a$, and $|W|>|U|$ 
	is infinite. This transfinite back and forth procedure
	gives us a pairing of elements (i.e., a bijection), 
	and as the entire $\dom_{0}(a)$ and $\dom_{i}(a)$ were enumerated, 
	$\dom(f_{0i}^a) = \dom_{0}(a)$, and $\ran(f_{0i}^a) = \dom_i(a)$. 

	If $a, b\subseteq\d_{0i}'$ are distinct atoms, then we have two cases 
	according to Lemma \ref{dom-particio}: either 
	\begin{align}
		\dom_0(a) = \dom_0(b),\quad\text{ and }\quad \dom_i(a) = \dom_i(b),
		\label{dom-case-1}
	\end{align}
	or
	\begin{align}
		\dom_0(a) \cap \dom_0(b) = \emptyset,\quad\text{ and } \quad
		\dom_i(a) \cap \dom_i(b) = \emptyset.
	\end{align}
	We construct $f_{0i}^a$ and $f_{0i}^b$ in such a way that
	in case \eqref{dom-case-1} we pick the same bijections: 
	If \eqref{dom-case-1} holds for the atoms $a$ and $b$, then
	we choose $f_{0i}^a = f_{0i}^b$.
	
	Define
	\begin{align}
		f_{0i} \defeq \bigcup\{ f_{0i}^a:\; a\in At,\ a\subseteq\d_{0i}'\},
		\quad
		f_{00}\defeq \id.
	\end{align}
	As the union of atoms below $\d_{0i}'$ is $\d_{0i}'$ itself 
	(the algebra is a union-complete BAO), 
	$f_{0i}$ is a relation with $\dom(f_{0i}) = \dom_{0}(\d_{0i}')$, and 
	$\ran(f_{0i}) = \dom_{i}(\d_{0i}')$.

	\begin{lemma}\label{lemma:n1}
		$f_{0i}$ is a bijective function from $V$ to $V$.
	\end{lemma}
	\begin{proof}[of Lemma \ref{lemma:n1}]
		Combining $\dom(f_{0i}) = \dom_{0}(\d_{0i}')$, and 
		$\ran(f_{0i}) = \dom_{i}(\d_{0i}')$ with Lemma \ref{egeszV}, we get
		$\dom(f_{0i}) = V$ and $\ran(f_{0i})=V$. 
		By Lemma \ref{dom-particio} and Lemma \ref{egeszV}, the sets
		\begin{align}
			\big\{ \dom(f_{0i}^a):\; a\in At, \ a\subseteq \d_{0i}'   \big\},\text{ and }
			\big\{ \ran(f_{0i}^a):\; a\in At, \ a\subseteq \d_{0i}'   \big\}
		\end{align}
		are partitions of $V$. Lemma \ref{dom-particio} also implies
		that it is not possible to have atoms $a, b\subseteq\d_{ij}'$
		such that $\dom_0(a)=\dom_0(b)$ but $\dom_i(a)\neq \dom_i(b)$; 
		or $\dom_0(a)\neq \dom_0(b)$ but $\dom_i(a) =\dom_i(b)$. Since in 
		the case $\dom_0(a)=\dom_0(b)$ and $\dom_i(a) = \dom_i(b)$ 
		we chose $f_{0i}^a = f_{0i}^b$, it follows that $f_{0i}$
		is the union of bijections having pairwise disjoint domains 
		and pairwise disjoint ranges. 
		Putting all this together, we get that $f_{0i}$ is a bijection. 
	\end{proof}
	
	In the next step, we extend $f_{0i}$ to a bijection
	$f: {}^{\alpha}V\to {}^{\alpha}V$ by letting for every
	$s\in {}^{\alpha}V$,
	\begin{align}
		f(s)\defeq \<  f_{0i}^{-1}\big(s(i)\big):\; i<\alpha \> \in {}^{\alpha}V\,.
	\end{align}
	This $f$ induces the mapping 
	\begin{align}
		G(X)\defeq f[X]=\{ f(x):x\in X\} 
	\end{align}
	
	\begin{lemma}\label{lemma:n0}
		The mapping $G$ is an isomorphism 
		\begin{align}
			G: \gA' = \< A', \cup, \setminus, \C_i^{V}, 
		\s_{ij}', \d_{ij}' \>_{i<\alpha}\quad \izomorf \quad
		\gA'' = \< A'', \cup, \setminus, \C_i^{V}, 
		\s_{ij}'', \d_{ij}''\>_{i<\alpha},
		\end{align}
		where $A''$, $\s_{ij}''$, and $\d_{ij}''$ are the $G$-images 
		of $\s_{ij}'$ and $\d_{ij}'$, respectively. Further, $G$ takes infinite unions to unions.
	\end{lemma}
	The content of the lemma is that the set theoretic operations 
	$\C_i^V$ are kept fixed under the isomorphism $G$. 
	This lemma essentially follows 
	from \cite[Theorem 5.1.34]{HMT}, nevertheless we provide the proof.\medskip
	\begin{proof}[of Lemma \ref{lemma:n0}]
		The case of Boolean operations is straightforward, and thus
		we only need to check $G(\C_i^V(R)) = \C_i^V(G(R))$ for 
		$R\subseteq {}^{\alpha}V$. Write $h=f^{-1}$. Then $G(X) = \{x: h(x)\in X\}$. 
		Unfolding the definitions yields
		\begin{align}
			s\in \C_i^V(G(R)) &\Leftrightarrow (\exists u\in V) h(s(i/u))\in R, \\
			s\in G(\C_i^V(R)) &\Leftrightarrow (\exists u\in V) h(s)(i/u)\in R.
		\end{align}
		Now, by definition of $f$ we have $h(s)(j) = f_{0j}(s(j))$, so
		\begin{align}
			h\big(s(i/u)\big)(j) = f_{0j}\big( s(i/u)(j) \big) = 
			 \begin{cases}
				f_{0j}\big( s(j) \big) & \text{ if } i\neq j, \\
				f_{0j}\big( u \big) & \text{ otherwise.}
			\end{cases}
		\end{align}
		On the other hand, 
		\begin{align}
			\big(h(s)(i/u)\big)(j)  = 
			 \begin{cases}
				f_{0j}\big( s(j) \big) & \text{ if } i\neq j, \\
				u & \text{ otherwise.}
			\end{cases}
		\end{align}
		As each $f_{0j}$ is a bijection, it follows that
		$(\exists u\in V) h(s(i/u))\in R$ holds if and only if 
		$(\exists u\in V) h(s)(i/u)\in R$ does.
		
		Finally, it is routine to check that $G$ takes infinite unions to unions.
	\end{proof}

	\begin{lemma}\label{lemma:n2}
		For every $i,j<\alpha$ we have $\D_{ij}^V\subseteq \d_{ij}''$.  (Note that possibly $\D_{ij}^V\notin A''$).
	\end{lemma}
	\begin{proof}[of Lemma \ref{lemma:n2}]
		We show first $\D_{0i}^V\subseteq \d_{0i}''$ for $i<\alpha$. Recall 
		from Lemma \ref{lemma:n1} that $\dom(f_{0i}) = \ran(f_{0i})=V$, 
		and by construction
		\begin{align}
		 	f_{0i}\times {}^{\alpha-2}V
			= \{s\in {}^{\alpha}V: s(0)\in\dom(f_{0i}),\ s(i)=f_{0i}(s(0))\}\ 
			 \subseteq\ \d_{0i}'\,.
			\label{eq:star}
		\end{align}
		(By cylindric axioms, $\C_k^V\d_{0i}' = \d_{0i}'$ for $k\neq 0, i$, 
		thus, sequences in $\d_{0i}'$ can take any values on the coordinates
		different from $0$ and $i$).
		Assume $s\in{}^{\alpha}V$, $s\in \D_{0i}^V$, that is, $s(0)=s(i)$. We need
		to prove $s\in \d_{0i}''$. As $f$ is a bijection, there is a unique
		$z\in {}^{\alpha}V$ such that $f(z)=s$. By definition of $f$, we
		have $s(i) = f_{0i}^{-1}\big( z(i)\big)$, in particular, $s(0)=z(0)$
		(because $f_{00}=\id$). Combining this with $s(0)=s(i)$, 
		we get $z(0) = f_{0i}^{-1}(z(i))$, that is,
		\begin{align}
			z(i) = f_{0i}(z(0))\,.
		\end{align}
		But this means that $z\in\d_{0i}'$, by \eqref{eq:star}. Therefore
		$s\in f(\d_{0i}') = \d_{0i}''$. We proved $\D_{0i}^V\subseteq \d_{0i}''$.
		
		The rest of the cases needs axioms of cylindric algebras.
		$\D_{i0}^V \subseteq \d_{i0}''$ because $\D_{i0}^V = \D_{0i}^V$, 
		$\D_{0i}^V\subseteq \d_{0i}''$, and $\CA_{\alpha}\models \d_{0i}=\d_{i0}$ (\cite[1.3.1]{HMT}), 
		in particular, $\d_{0i}'' = \d_{i0}''$. Finally, let $i,j<\alpha$, $i,j,0$ distinct. Then
		\begin{align}
			\D_{ij}^V = \C_0^V(\D_{i0}^V\cap \D_{0j}^V) \subseteq
			\C_0^V(\d_{i0}''\cap \d_{0j}'') = \d_{ij}''
		\end{align}
		by cylindric axiom (C6). 
	\end{proof}

	\medskip	
	\noindent So far, we have a cylindric algebra
	\begin{align}
		\gA'' = \< A'', \cup, \setminus, \C_i^{V}, 
		\s_{ij}'', \d_{ij}''\>_{i<\alpha},
	\end{align}
	such that $\Rd_{cs}(\gA'')\cong\Rd_{cs}(\gA)$; and $\gA''$ is
	complete, atomic, suprema are unions, and each $\s_{ij}''$ is 
	completely additive satisfying
	\begin{align}
		\s_{ij}''(x) = \begin{cases}
			x & \mbox{ if } i=j,\\
			\C_i^V(x\cap \d_{ij}'') & \mbox{ if } i\neq j,
		\end{cases} \label{eq:s=cd2}
	\end{align}
	and $\D_{ij}^V\subseteq \d_{ij}''$.
	\medskip
	
	Next we prove that $\s_{ij}''$ are the real, set-theoretic 
	operations $\S_{ij}^V$.

	\begin{lemma}\label{lemma:n3}
		For all $i,j<\alpha$, \ $\s_{ij}'' = \S_{ij}^V$.
	\end{lemma}
	\begin{proof}[of Lemma \ref{lemma:n3}]
		The case $i=j$ is straightforward: $\s_{ii}'' = \id = \S_{ii}^V$.
		Take $i\neq j<\alpha$. 		
		As the algebra is complete, atomic, and $\s_{ij}''$ is completely 
		additive, $\s_{ij}''$ is 
		determined by its values on atoms. Recall the definitions
		\begin{align}
			\C_i^V(a) &\defeq 
				\{ x\in {}^{\alpha}V:\; x(i/v)\in a\text{ for some }
				 v\in V\},\\
			 \S_{ij}^V(a) &\defeq 
			 	\{ x\in {}^{\alpha}V:\; x(i/x(j))\in a \}.
		\end{align}
		By \eqref{eq:s=cd2} we have $\s_{ij}''(x) = \C_i^V(x\cap \d_{ij}'')$.
		
		For an atom $a\not\subseteq\d_{ij}$ we have $\s_{ij}''(a) = \emptyset$ by 
		$a\cap\d_{ij}=\emptyset$, \eqref{eq:s=cd2} and $\C_i^V\emptyset=\emptyset$. 
		As $\D_{ij}^V\subseteq\d_{ij}''$ we also get $\S_{ij}^V(a)=\emptyset$. 
		Therefore, what we need to check is whether
		\begin{align}
			\s_{ij}''(a) = 
			\S_{ij}^V(a)\quad\text{ for } a\subseteq \d_{ij}''\text{ atoms}.
		\end{align}
		Now, if $a$ is an atom below $\d_{ij}''$, 
		then $\s_{ij}''(a) = \C_i^V(a\cap\d_{ij}'') = \C_i^V(a)$, 
		and so we need to prove
		\begin{align}
			\C_i^V(a) =	\S_{ij}^V(a)\quad
			\text{ for } a\subseteq \d_{ij}''\text{ atoms}.
		\end{align}
		The inclusion $\S_{ij}^V\subseteq \C_i^V$ is immediate from the 
		definitions. For the other direction, take an arbitrary 
		$s\in \C_i^V(a)$, and let
		\begin{align}
			z \defeq s\big( i/s(j) \big)\,.
		\end{align}
		To get $s\in \S_{ij}^V(a)$ we need $z\in a$. Now, $z\in D_{ij}^V$
		as $z(i)=s(j)=z(j)$. As $\D_{ij}^V\subseteq \d_{ij}''$, we get
		$z\in \d_{ij}''$.
		
		Suppose, by way of contradiction, that $z\notin a$. Then
		$z\in -a\cap \d_{ij}''$, and thus $s\in \C_i^V(-a\cap \d_{ij}'')$.
		On the other hand, $s\in \C_i^V(a)$, and $a\subseteq \d_{ij}''$, 
		so $s\in \C_{i}^V(a\cap \d_{ij}'')$. It follows that
		\begin{align}
			s\in \C_i^V(-a\cap\d_{ij}'')\cap 
			\C_i^V(a\cap\d_{ij}''). 
		\end{align}
		But this contradicts cylindric axiom (C7). 
		Consequently, $z\in a$, and this completes the proof of Lemma~\ref{lemma:n3}.
	\end{proof}
	
	\noindent Summing up, we have shown that
	\begin{align}
		\< A'', \cup, \setminus, \C_i^{V}, 	\s_{ij}''\>_{i<\alpha} = 
		\< A'', \cup, \setminus, \C_i^{V}, 	\S_{ij}^V\>_{i<\alpha} \in
		\RA^{cs}_{\alpha}.
	\end{align}
	This implies $\gA\in\RA^{cs}_{\alpha}$, which makes the proof of Theorem
	\ref{thm:1a} complete.
\end{proof}

As a corollary to Theorem \ref{thm:1a}, we can fill in the status of the open line from $c$ to $cs$ 
in \cite[Figure 1]{AAPALI97ec}: it should be a dashed line indicating finite axiomatizability. 

\begin{theorem}\label{thm:2}  
	$\RA_{\alpha}^{cs}$ is finitely axiomatizable over $\RA_{\alpha}^{c}$ (by the  equations \ref{F:0}-\ref{F:last},(FS) defining $\FSC_{\alpha}$).
\end{theorem}
\begin{proof} This is an immediate corollary of Lemma~\ref{lem:1} and $\RA_{\alpha}^{cs}\subseteq\FSC_{\alpha}$.
\end{proof}

We now prove Theorem~\ref{main2-thm}. Recall first that $\RDF_{\alpha}=\mathbf{SP}\RA_{\alpha}^c$ and $\FRPA_{\alpha}=\mathbf{SP}\RA_{\alpha}^{csp}$. 
Then $\SRd_{cs}\FRPA_{\alpha}=\mathbf{S}\mathbf{SP}\Rd_{cs}\RA_{\alpha}^{csp}=\mathbf{SP}\RA_{\alpha}^{cs}$. So Theorem~\ref{main2-thm} states that $\gA\in\mathbf{SP}\RA_{\alpha}^{cs}$ iff ($\Rd_c\gA\in\mathbf{SP}\RA_{\alpha}^c$  and equations \eqref{main6}--\eqref{main9} hold in $\gA$). By Theorem~\ref{thm:2} we get that $\gA\in\RA_{\alpha}^{cs}$ iff ($\gA\in\RA_{\alpha}^c$  and equations \ref{F:0}-\ref{F:last},(FS) hold in $\gA$).  Equations \eqref{main6}--\eqref{main9} in Theorem \ref{main2-thm} are the equations in \ref{F:0}-\ref{F:last},(FS) in which the substitutions occur%
\footnote{The indices $i,j,k$ in Theorem \ref{main2-thm} are required to be distinct while in (F0)--(F6),(FS) they are not required to be distinct. It is easy to check that this does not make any difference.}; 
the ones in which the substitutions do not occur hold in $\RA_{\alpha}^c$, so they are superfluous. Theorem~\ref{main2-thm} has been proved.

\section{$\RA_{\alpha}^{csp}$ is finitely axiomatizable over $\RA^{cs}_{\alpha}$}
\label{sec:povercs}

This section proves our second main result which states that a finitary polyadic algebra of finite dimension is representable iff its permutation-free reduct is representable. A corollary to this 
result is that permutations over the Booleans together with substitutions and cylindrifications are finitely axiomatizable (by the axioms of finitary polyadic algebras). 

\begin{theorem}\label{thm:2a}
		For finite $\alpha>2$, an $\alpha$-dimensional algebra 
		$\gA\in\FPA_{\alpha}$ is 
		representable if and only if 
		its permutation-free reduct is representable.
\end{theorem}
\begin{proof}
	The left-to-right direction is immediate. We prove the converse direction.
	The beginning of the proof follows the proof of Theorem \ref{thm:1a}.
	Namely, by the proof of Lemma \ref{lem:1}, with the obvious modifications,
	it is enough to prove
	\begin{align*}
		&\text{if } \gA\in\FPA_{\alpha} \text{ is union-complete, and } 
		\Rd_{cs}\gA = \<A,\cup,\setminus,\C_i^U, \S_{ij}^U\>_{i,j<\alpha},\\
		&\text{ for some } U\text{ and } 
		A\subseteq\Sb({}^{\alpha}U), \text{ then } \gA\in\RA^{csp}_{\alpha}
		\tag{$\sharp$'}
		\label{uj:*}
	\end{align*}
	(In showing $g(x)=0\Rightarrow g(\p_{ij}x)=0$ one uses (F6) in place of (F3).)\\
	
	To prove \eqref{uj:*}, take a union-complete $\gA\in\FPA_{\alpha}$
	such that \[\Rd_{cs}\gA = \<A,\cup,\setminus,\C_i^U, \S_{ij}^U\>_{i,j<\alpha}\]
	for some base set $U$ and $A\subseteq \Sb({}^{\alpha}U)$.
	Putting back the permutation operations, we get
	\begin{align}
		\gA = \< A, \cup, \setminus, \C_i^{U}, \S_{ij}^{U},
			\p_{ij}^{\gA} \>_{i,j<\alpha},
	\end{align}
	where the operations $\p_{ij}^{\gA}$ are abstract, not yet set-represented 
	operations that we want to represent by the operations $\P_{ij}$ over sets.

As in the analogous step in the proof of Theorem \ref{thm:1a}, we can find elements $\d_{ij}^\gA\in A$ for $i,j<\alpha$ such that
		$\<A, \cup, \setminus, \C_i^U, \d_{ij}^{\gA}\>_{i,j<\alpha}$ is a 
	cylindric algebra, and
	\begin{align}
		\S_{ij}^{U}x = \begin{cases}
			x & \mbox{ if } i=j,\\
			\C_i^U(x\cap \d_{ij}^{\gA}) & \mbox{ if } i\neq j
		\end{cases} \label{uj:eq:s=cd}
	\end{align}
	holds. We also have
	\begin{align}
		\D_{ij}^U\subseteq\d_{ij}^{\gA} \label{eq:Dij:dij}
	\end{align}
	since \eqref{uj:eq:s=cd} implies $\S_{ij}^U(-\d_{ij}^{\gA})=\emptyset$, so $\D_{ij}^U\cap -\d_{ij}^{\gA}=\emptyset$ by the definition of $\S_{ij}^U$.
	
	\begin{lemma}\label{lem:pijdkl}
		$\p_{ij}^{\gA}(\d_{kl}^{\gA}) = \d_{\tau(k)\tau(l)}^{\gA}$ 
		where $\tau = [i,j]$. 
	\end{lemma}
	\begin{proof}[of Lemma \ref{lem:pijdkl}]
		For any $i,j<\alpha$, 
		from 
		\eqref{uj:eq:s=cd}, \cite[1.5.4(i)]{HMT} and cylindric axiom (C2) we get
		$\S_{ij}^U(\d_{ij}^{\gA}) = 1$ and
		$x\cap\d_{ij}^\gA \subseteq \S_{ij}^U(x)$. Sain--Thompson \cite{Sain-Tho}
		define the class $\FPEA_{\alpha}$ as $\FPA_{\alpha}$ enriched
		with diagonal elements $\d_{ij}$ satisfying $\s_{ij}\d_{ij}=1$ and 
		$x\cdot \d_{ij}\leq \s_{ij}x$ (see \cite[p.544, Definition 2]{Sain-Tho}). 
		We just verified that these two properties
		hold for $\S_{ij}^U$ and $\d_{ij}^{\gA}$, therefore the expanded 
		structure $\<\gA, \d_{ij}^{\gA}\>_{i,j<\alpha}$ is an $\FPEA_{\alpha}$.
		\cite{Sain-Tho} proves that $\FPEA_{\alpha}$ are term equivalent 
		to quasi-polyadic equality algebras $\QPEA_{\alpha}$, and in particular, 
		\cite[p.557-558, (E3)]{Sain-Tho} proves that from $\FPEA_{\alpha}$ the identity
		\begin{align}
			\p_{ij}(\d_{kl}) = \d_{\tau(k)\tau(l)}
		\end{align}
		follows.
	\end{proof}

	Let 
	$\gB$ be the Boolean algebra generated by the elements $\C_i^U(x)$ for $i<\alpha$
	and $x\in A$:
	\begin{align}
		\gB \defeq \Sg^{\Bl\gA}\big( \{ \C_i^U(x):\; x\in A,\  i<\alpha \}  \big)\,.
	\end{align}

	\begin{lemma}\label{uj:pij-igazi-Bn}
		$\p_{ij}^{\gA}(x) = \P_{ij}^U(x)$ for every $x\in B$ and $i,j<\alpha$.
	\end{lemma}
	\begin{proof}[of Lemma \ref{uj:pij-igazi-Bn}]
		For $i=j$ the statement follows from \ref{FPA:0} $\p_{ii}^{\gA}=\id=\P_{ii}^U$.
		Let $i\neq j$. 
		As $\p_{ij}^{\gA}$ is a Boolean endomorphism by \ref{FPA:endo}, it is enough to
		check 
		\begin{align}
			\p_{ij}^{\gA}(\C_k^U(a)) = \P_{ij}^U(\C_k^U(a))
		\end{align}
		for $a\in A$ and $k<\alpha$. We have three cases depending on whether
		$k=i$, $k=j$, or $k\neq i,j$. The
		identities \ref{uj:p2}, \ref{uj:p3}, and \ref{uj:p1} can
		be used for these cases. For example, if $k=i$, then
		\begin{align}
			\p_{ij}^{\gA}\big(\C_i^U(a)\big) 
			\overset{\text{\ref{uj:p2}}}{=} \S_{ji}^U(\C_i^U(a)) = \P_{ij}^U(\C_i^U(a)),
		\end{align}
		where the last equality is a property of the set theoretic $\P_{ij}^U$.
		The rest of the cases is similar.
	\end{proof}

	\begin{lemma}\label{uj:pij-igazi-dij-alatt}
		$\p_{ij}^{\gA}(a) = \P_{ij}^U(a)$ for every $a\in A$, $a\subseteq \d^{\gA}_{kl}$, $k\neq l<\alpha$.
	\end{lemma}
	\begin{proof}
		Recall that $\<A, \cup, \setminus, \C_i^U, \d_{ij}^{\gA}\>_{i,j<\alpha}$ is a cylindric algebra.
		Choose $m<\alpha$ such that $m\notin\{k,l\}$; this is possible since $\alpha\geq 3$. 
		By \cite[1.3.3]{HMT} we have
		\[
			\C_m^{U}\d_{kl}^{\gA}=\d_{kl}^{\gA},
		\]
		thus the definition of $\gB$ yields $\d_{kl}^{\gA}\in B$. 
		Moreover, since $a\subseteq \d_{kl}^{\gA}$, the cylindric algebra identity
		\[
			\C_k^{U}(a)\cap \d_{kl}^{\gA} = \C_k^U(a\cap \d_{kl}^{\gA})\cap \d_{kl}^{\gA} = a\cap \d_{kl}^{\gA} = a
		\]
		holds (the penultimate equation follows from \cite[1.3.9]{HMT}). 
		Now $\C_k^{U}(a)\in B$ by the definition of $B$, and $\d_{kl}^{\gA}\in B$. As $\gB$ is a Boolean subalgebra of $\gA$, 
		it follows that $a\in B$. Hence, by Lemma \ref{uj:pij-igazi-Bn},
		\[
			\p_{ij}^{\gA}(a)=\P_{ij}^{U}(a),
		\]
		as required.
\end{proof}

	\medskip

	Recall (see the paragraph before Lemma \ref{lem:permulem}) 
	that to every permutation $\sigma:\alpha\to\alpha$ 
	there is a corresponding derived polyadic operation $\p_{\sigma}$
	such that $\p_{[i,j]} \defeq \p_{ij}$ and $\p_{\sigma\circ \tau} = \p_{\sigma}\circ \p_{\tau}$. 
	The set-theoretic permutation operator is defined as
	\[
		\P_{\sigma}(x) = \big\{ s\circ\sigma^{-1}:\; s\in x \big\} .
	\]
	It is easy to see that Lemmas \ref{lem:pijdkl}, \ref{uj:pij-igazi-Bn} and 
	\ref{uj:pij-igazi-dij-alatt} extend for any permutation $\sigma$ in the
	following manner:
	\begin{align}
		\p_{\sigma}^{\gA}(\d_{kl}^{\gA}) &= \d_{\sigma(k)\sigma(l)}^{\gA}, \label{eq:pdd} \\
		\p_{\sigma}^{\gA}(x) &= \P_{\sigma}^U(x)\quad 
			\text{ for every } x\in B\,,  \label{eq:aha}   \\
		\p_{\sigma}^{\gA}(a) &= \P_{\sigma}^U(a)\quad
			 \text{ for every } a\in A, \ 
			a\subseteq \d_{kl}^{\gA},\ k\neq l<\alpha\,. \label{eq:pPad}
	\end{align}
	Let $PT_{\alpha}\subset T_{\alpha}$	denote the set of permutations of $\alpha$. 
	Write $At$ for the atoms of $\gA$.
	\medskip

	Call an atom $a\in A$ \emph{repetition-free} if $a\not\subseteq \d_{kl}^{\gA}$ 
	for every $k\neq l<\alpha$.  
	Write
	\begin{align}
		\RfAt = \big\{ a\in At :\; (\forall k\neq l<\alpha)\ 
		a\not\subseteq \d_{kl}^{\gA}\   \big\}\,
	\end{align}
	and let 
	$S$ denote the set of all sequences belonging to repetition-free atoms
	\begin{align}
		S &= \bigcup \big\{ a:\; a\in\RfAt \big\}\,.
	\end{align}
	A sequence $s\in {}^{\alpha}U$ is called repetition-free if $s\notin \D_{kl}^{U}$ for every $k\neq l<\alpha$. 
	Elements of $S$ are repetition-free sequences by \eqref{eq:Dij:dij}, but $S$ is 
	not necessarily the set of all such sequences, because atoms below some
	$\d_{kl}^{\gA}$ might also contain repetition-free sequences. 
	We also have that $S$ is closed under permutations, i.e., for $s\in S$ and $\sigma\in PT_{\alpha}$, 
	we have $s\circ\sigma\in S$.  This is so, because $s\circ\sigma\notin S$ means that $s\circ\sigma\in a\subseteq \d_{kl}^{\gA}$ for some $k,l\in\alpha$. Then $s\in\P_{\sigma}(a) = \p_{\sigma}(a) \subseteq \d_{\sigma(k)\sigma(l)}^{\gA}$ by \eqref{eq:pPad} and \eqref{eq:pdd}, so $s\notin S$ by the definition of $S$.
	For each $s\in S$ define its block as 
	\begin{align}
		PT(s) = \{ s\circ\sigma:\; \sigma\in PT_{\alpha} \}
	\end{align}
	and let $R$ contain a representant element from each block, 
	that is $|R\cap PT(s)|=1$ for each $s\in S$. 
	Then 
	\begin{align}
		S = \{ s\circ\sigma:\; s\in R\mbox{ and }\sigma\in PT_{\alpha} \},
		\end{align}
	and for each $z\in S$ there are unique $s\in R$ and $\sigma\in PT_{\alpha}$ such that $z=s\circ\sigma$.
	\bigskip
	
Our task is to represent $\p_{\sigma}$ on the atoms by $\P_{\sigma}$. By \eqref{eq:pPad}, we have to change the representation only on the repetition-free atoms. 
	How can it happen that $\p_{\sigma}(a)\ne\P_{\sigma}a$?
	For example, it may happen that $s\circ\sigma\in a\in \RfAt$, then we would need $s\in\p_{\sigma}a$, but $s\in b\ne \p_{\sigma}a$. Instead of ``pulling $s$ out of $b$ and putting it into $p_{\sigma}a$", we make a ``copy" $s_{\sigma}$ of $s$ and put $s_{\sigma}$ into $\p_{\sigma}(a)$. In more detail, we will make copies $s_{\sigma}$ of $s\in R$ for all $\sigma\in PT_{\alpha}$ and put $s_{\sigma}\circ\delta^{-1}$ to $\p_{\delta}\p_{\sigma}(a)$ whenever $s\circ\sigma\in a$. Then, as we will show, $\p_{\sigma}$ will be well-represented on $\RfAt$ and the representations of the atoms will remain disjoint. How do we make copies $s_{\sigma}$ of $s$? First, we change $s$ into a set $\hat{s}$ of sequences as in the previous section by replacing $U$ with $U\times W$ for a suitable $W$, and then we will partition the set $\hat{s}$ to obtain the $s_{\sigma}$. We now begin to execute this plan.

	We change the representation by replacing $s\in {}^\alpha U$ with $\hat{s}\subseteq {}^{\alpha}(U\times W)$ as in the previous section, see \eqref{eq:hats}. 
	The new base set of the representation will be $V\defeq U\times W$ where $W=PT_{\alpha}$. 
	Then ${}^{\alpha}V=\{ s\times w : s\in{}^{\alpha}U \mbox{ and } w\in{}^{\alpha}W\}$.  For each sequence $s\in {}^{\alpha}U$ let  
	\begin{align}
		\hat{s} \defeq \{ s\times w :\;  w\in {}^{\alpha}W \},\label{eq:shatdef}
	\end{align}
	and for $X\in\Sb({}^{\alpha}U)$ let
	\begin{align}
		F(X) \defeq \bigcup\big\{ \hat{s}:\; s\in X \big\}\in\Sb({}^{\alpha}V). \label{def:F}
	\end{align}
	It is straightforward to check that $F$ is an embedding 
	\begin{align}
		\< \Sb({}^{\alpha}U), \cup, \setminus, \C_i^{U}, \S_{ij}^U,  \P_{ij}^U\,\>_{i<\alpha}
		 \quad\rightarrowtail\quad
		\< \Sb({}^{\alpha}V), \cup, \setminus, \C_i^{V}, \S_{ij}^V, \P_{ij}^V\, \>_{i<\alpha}.
	\end{align}
	Consequently, there is a subset $A'\subseteq \Sb({}^{\alpha}V)$ 
	such that $\<\gA, \d_{ij}^{\gA}\>_{i,j<\alpha}$ is isomorphic to 
	an algebra
	\begin{align}
		\gA' = 
		\< A', \cup, \setminus, \C_i^{V}, 
		\S_{ij}^V, \p_{ij}', \d_{ij}' \>_{i,j<\alpha},
	\end{align}
	where $\p_{ij}'$ and $\d_{ij}'$ are the isomorphic 
	images of $\p_{ij}^{\gA}$ and $\d_{ij}^{\gA}$, respectively. 
	Let $At'$ denote the set of atoms of $\gA'$, that is, the embedded image
	of $At$. Similarly, $\gB'$ is the embedded image of $\gB$. 
	For every permutation $\sigma$, equations  
	\eqref{eq:pdd}-\eqref{eq:pPad} ensure
	\begin{align}
		\p_{\sigma}'(\d_{kl}') &= \d_{\sigma(k)\sigma(l)}' \label{eq:pdd2} \\
		\p_{\sigma}'(x) &= \P_{\sigma}^V(x)\quad 
			\text{ for every } x\in B'\,,\\
		\p_{\sigma}'(a) &= \P_{\sigma}^V(a)\quad
			 \text{ for every } a\in A', \ 
			a\subseteq \d_{kl}',\ k\neq l<\alpha\,.\label{eq:repatomokonjo}
	\end{align}

	\noindent
	Next, for each $s\in R$ we split the embedded image $\hat{s}$ ($=F(\{s\})$)
	into $\alpha!$ ($\alpha$ factorial) elements.
	
	\begin{lemma}\label{huha-lemma} 
		There is a partition $\{ s_{\sigma} : \sigma\in PT_{\alpha}\}$ of $\hat{s}$ such that
		$\C_i^V(s_{\sigma}) = \C_{i}^V(\hat{s})$ for all $\sigma\in PT_{\alpha}$ 
		and $i<\alpha$ .
		\end{lemma}
	\begin{proof} Recall that $W=PT_{\alpha}$. Let $\star$ be any operation on $W$ such that $\langle W,\star\rangle$ is a 
		group, e.g., we can take $\star$ to be a bijective image of addition modulo $|W|$ on $\{ 0,1,\dots,|W|-1\}$.  For $\sigma\in W$ define
		\begin{align}
			s_{\sigma}\defeq\{ s\times w\in{}^\alpha V : w_0\star\dots\star w_{\alpha-1}=\sigma\}.
			\end{align}
	It is now routine to check that $\{ s_{\sigma} : \sigma\in PT_{\alpha}\}$ satisfies the  requirements.
	\end{proof}
	\medskip

	Let us change the representation $F(a)$ of the repetition-free atoms $a$ to $rep(a)$ defined as
	\begin{align}
		\rep(a)\defeq \bigcup\big\{ \P^V_{\delta^{-1}}(s_{\sigma}):\; 
		\sigma,\delta\in PT_{\alpha},\ \ 
		 s\circ\sigma\in\p^{\gA}_{\sigma^{-1}\circ\delta}(a)\mbox{ and }s\in R  
		  \big\} .\label{eq:rep}
	\end{align}
	\medskip
	
	\begin{lemma}\label{lem:P-re-jo}
		For $a\in \RfAt$ and $\tau\in PT_{\alpha}$ we have
		\begin{align}
			\rep\big( \p_{\tau}^{\gA}(a) \big) = \P_{\tau}^V(\rep(a))\,.
		\end{align}
	\end{lemma}
	\begin{proof}[of Lemma \ref{lem:P-re-jo}]
		\begin{align}
			\rep(\p_{\tau}^{\gA}(a)) &= \bigcup\big\{
				 \P^V_{\delta^{-1}}(s_{\sigma}):\; 
				s\circ\sigma\in\p^{\gA}_{\sigma^{-1}\circ\delta}(\p_{\tau}^{\gA}a),\ s\in R
		  	  	\big\} \\
			&= \bigcup\big\{ \P^V_{\delta^{-1}}(s_{\sigma}):\; 
				s\circ\sigma\in
				\p^{\gA}_{\sigma^{-1}\circ\delta\circ\tau}(a),\ s\in R  \big\} \\
			&= \bigcup\big\{ \P^V_{\tau}\P^V_{\tau^{-1}}\P^V_{\delta^{-1}}(s_{\sigma}):\; 
				s\circ\sigma\in
				\p^{\gA}_{\sigma^{-1}\circ\delta\circ\tau}(a),\ s\in R  \big\} \\
			&= \P^V_{\tau}\bigcup\big\{\P^V_{(\delta\circ\tau)^{-1}}(s_{\sigma}):\; 
				s\circ\sigma\in
				\p^{\gA}_{\sigma^{-1}\circ\delta\circ\tau}(a),\ s\in R  \big\} \\
			&= \P^V_{\tau}(\rep(a)) .
		\end{align}
	\end{proof}

Recall that $F$ is defined in \eqref{def:F}. 
	\begin{lemma}\label{lem:C-re-jo}
		For $a\in \RfAt$ and $i<\alpha$ we have
		\begin{align}
			\C_i^V(\rep(a)) = \C_i^V(F(a))\,.\label{eq:uccso}
		\end{align}		
	\end{lemma}
	\begin{proof}[of Lemma \ref{lem:C-re-jo}]
		($\supseteq$) Suppose $z\in \C_i^V(F(a))$ for some $a\in \RfAt$ and
		$i<\alpha$. There are a repetition-free sequence $s\in R$ and a
		permutation $\sigma$ such that $s\circ \sigma\in a$, and 
		$z\in \C_i^V(F(\{s\circ\sigma\}))$. With $\delta=\sigma$ in the
		definition \eqref{eq:rep} of $\rep(a)$, and using $s\circ\sigma\in a$
		we get that $\P^V_{\sigma^{-1}}(s_{\sigma})\subseteq \rep(a)$. Observe
		that
		\begin{align}
			F(\{s\circ \sigma\}) = \P^V_{\sigma^{-1}}(F(\{s\}))\label{eq:fpfs}
		\end{align}
		because $F$ is an embedding with respect to the set-theoretic 
		polyadic operations. By construction, $s_{\sigma}\subseteq F(\{s\})$, 
		and $\C_j^V(s_{\sigma}) = \C_j^V(F(\{s\}))$ for all $j<\alpha$. Therefore, 
		\begin{eqnarray}
		  z\in \C_i^V(F(\{s\circ\sigma\})) 
			&\overset{\text{\eqref{eq:fpfs}}}{=}& 
				\C_i^V(\P^V_{\sigma^{-1}}(F(\{s\}))) \label{eq:x1}\\
			&\overset{\text{\ref{lemma1.6-kov}}}{=}&
				 \P^V_{\sigma^{-1}}\C_{\sigma^{-1}(i)}^V(F(\{s\}))\\
			&=& \P^V_{\sigma^{-1}}\C_{\sigma^{-1}(i)}^V( s_{\sigma} )\\
			&\overset{\text{\ref{lemma1.6-kov}}}{=}& 
				\C_i^V(\P^V_{\sigma^{-1}}( s_{\sigma} )) \label{eq:x2}\\
			&\subseteq& \C_i^V(\rep(a)).		
		\end{eqnarray}
		\medskip
		
		\noindent ($\subseteq$) In the other direction we need that
		$\C_i^V(\rep(a)) \subseteq \C_i^V(F(a))$. By definition, and 
		complete additivity of $\C_{i}^V$, 
		\begin{align}
			\C_i^V(\rep(a)) = \bigcup\big\{ \C_{i}^V( \P^V_{\delta^{-1}}(s_{\sigma})):
				\; s\circ\sigma\in 
				\p^{\gA}_{\sigma^{-1}\circ\delta}(a),\ s\in R
			\big\} .
		\end{align}
		Take $s\in R$, $\sigma, \delta\in PT_{\alpha}$ such that 
		$s\circ\sigma\in \p^{\gA}_{\sigma^{-1}}\p^{\gA}_{\delta}(a)$.
		Now, \begin{align}
					\C_{i}^V(\P^V_{\delta^{-1}}(s_{\sigma})) = 
					\C_{i}^V(\P^V_{\delta^{-1}}(F(\{s\})))
				\end{align}
		(cf.\ steps \eqref{eq:x1}-\eqref{eq:x2}). As $F$
		is a homomorphism with respect to the ``real'' permutations, 
		\begin{align}
			\P^V_{\delta^{-1}}(F(\{s\})) = F(\{ s\circ \delta\}),
		\end{align}
		and thus we get $\C_i^V(\P^V_{\delta^{-1}}(s_{\sigma})) = \C_i^V(F(\{s\circ\delta\}))$.

		To complete the proof, it is enough to show that 
		$s\circ \sigma\in \p^{\gA}_{\sigma^{-1}}\p^{\gA}_{\delta}(a)$
		implies $\C_i^V(F(\{s\circ\delta\}))\subseteq \C_i^V(F(a))$. 
		Let $\eta = \delta^{-1}\circ\sigma$, and
		$b = \p^{\gA}_{\sigma^{-1}}\p^{\gA}_{\delta}(a)$.
		Now, 
		\begin{eqnarray}
			s\circ\sigma \in b &\Longrightarrow& s\in \P_{\sigma}(b) \\
			 &\Longrightarrow& s\circ\delta\in \P_{\delta^{-1}}\P_{\sigma}(b)
			  = \P_{\eta}(b).
		\end{eqnarray}
		Then
		\begin{align}
			&\C_i^V(F(\{s\circ\delta\}))\subseteq \C_i^V\P^V_{\eta}(F(b))
			\overset{\text{\ref{lemma1.6-kov}}}{=} \P^V_{\eta}\C_{\eta(i)}^V(F(b)) 
			= F(\P^U_{\eta}\C_{\eta(i)}^U(b))
			\\  & \overset{\text{\eqref{eq:aha}}}{=} 
			F(\p_{\eta}^{\gA}\C_{\eta(i)}^U(b)) 
			\overset{\text{\ref{lemma1.6-kov}}}{=}
			F(\C_{i}^U\p_{\eta}^{\gA}(b)) = F(\C_i^U(a)) = \C_i^V(F(a))\,.
		\end{align}
	\end{proof}

We are ready to define our final representation 
\[ G:\langle A,\cup,\setminus,\C_i^U,\S_{ij}^U, \p_{ij}^{\gA}\rangle_{i,j<\alpha}\to \langle\Sb({}^{\alpha}V),\cup,\setminus,\C_i^V,\S_{ij}^V, \P_{ij}^V\rangle_{i,j<\alpha} .\]
\noindent For atoms $a\in At$ we define 
	\begin{align}
		G(a) = \begin{cases}
			F(a) & \text{ if } a\subseteq \d_{kl}\text{ for some } k\neq l<\alpha,\\
			\rep(a) & \text{ if } a\in \RfAt\,.
		\end{cases}
	\end{align}
	$G$ maps atoms of $\gA$ into elements of $\Sb({}^{\alpha}V)$.

	\begin{lemma}\label{lem:G-Boolean}
		The sets $\{G(a):a\in At\}$
		are nonempty, pairwise disjoint, and their union is ${}^{\alpha}V$.
		Consequently, the map $G$ defined on the atoms extends to a Boolean
		embedding
		\begin{align}
			\< A,\cup,\setminus\>\quad\rightarrowtail\quad
			\<\Sb({}^{\alpha}V),\cup,\setminus\>
		\end{align}
		by
		\begin{align}
			G(x)=\bigcup\big\{G(a):a\leq x,\ a\in At\big\}.
			\label{eq:G-extension}
		\end{align}
	\end{lemma}
	\begin{proof}
		Since $\gA$ is complete and atomic and its suprema are set-theoretic
		unions, the atoms in $At$ form a partition of ${}^{\alpha}U$.
		In particular, every sequence in ${}^{\alpha}U$ belongs to a unique atom.

	We first show that the sets $\rep(a)$, for $a\in\RfAt$, form a partition
	of $F(S)$. For $s\in R$ and $\sigma,\delta\in PT_\alpha$, put
	\begin{align}
		E(s,\sigma,\delta)
		\defeq \P^V_{\delta^{-1}}(s_\sigma).
		\label{eq:E-cell}
	\end{align}
	Fix $s\in R$ and $\delta\in PT_\alpha$. By Lemma~\ref{huha-lemma},
	$\{s_\sigma:\sigma\in PT_\alpha\}$ is a partition of $\hat{s}$.
	Since $\P^V_{\delta^{-1}}$ is a bijection on ${}^{\alpha}V$,
	\[
		\{E(s,\sigma,\delta):\sigma\in PT_\alpha\}
	\]
	is a partition of
	\begin{align}
		\P^V_{\delta^{-1}}(\hat{s})
		&=\P^V_{\delta^{-1}}(F(\{s\}))
		  =F(\{s\circ\delta\}).
		\label{eq:E-fiber}
	\end{align}
	Recall that, by the choice of $R$, for every $z\in S$ there are unique
	$s\in R$ and $\delta\in PT_\alpha$ such that $z=s\circ\delta$.
	Hence, as $s$ and $\delta$ vary, the sets
	$F(\{s\circ\delta\})$ form a partition of $F(S)$. It follows from
	\eqref{eq:E-fiber} that
	\begin{align}
		\big\{E(s,\sigma,\delta):
		s\in R,\ \sigma,\delta\in PT_\alpha\big\}
		\label{eq:E-partition}
	\end{align}
	is a partition of $F(S)$.

	We next show that every $E(s,\sigma,\delta)$ in \eqref{eq:E-partition} occurs in
	$\rep(a)$ for exactly one $a\in\RfAt$. Fix
	$s\in R$ and $\sigma,\delta\in PT_\alpha$, and put $\eta=\sigma^{-1}\circ\delta$.
	Let $b$ be the unique atom containing $s\circ\sigma$.
	Since $s\in S$ and $S$ is closed under permutations,
	$s\circ\sigma\in S$, so $b\in\RfAt$.

	The operation $\p_\eta^{\gA}$ is a Boolean automorphism and hence
	permutes the atoms of $\gA$. Therefore there is a unique atom $a$ such
	that
	\begin{align}
		\p_\eta^{\gA}(a)=b.
		\label{eq:unique-a}
	\end{align}
	This atom $a$ is repetition-free. Indeed, if
	$a\subseteq\d_{kl}^{\gA}$ for some $k\neq l$, then by
	\eqref{eq:pdd}, we have
	\[
		b=\p_\eta^{\gA}(a)
		\subseteq
		\p_\eta^{\gA}(\d_{kl}^{\gA})
		=\d_{\eta(k)\eta(l)}^{\gA},
	\]
	contrary to $b\in\RfAt$.

	By \eqref{eq:unique-a}, we have
	\[
		s\circ\sigma\in
		\p_{\sigma^{-1}\circ\delta}^{\gA}(a),
	\]
	and hence, by the definition \eqref{eq:rep} of $\rep$,
	\begin{align}
		E(s,\sigma,\delta)\subseteq\rep(a).
		\label{eq:E-assigned}
	\end{align}
	Moreover, $a$ is the only atom for which
	\eqref{eq:E-assigned} holds by way of this occurrence in
	\eqref{eq:rep}: if
	\[
		s\circ\sigma\in
		\p_{\sigma^{-1}\circ\delta}^{\gA}(a')
	\]
	for another atom $a'$, then both
	$\p_{\sigma^{-1}\circ\delta}^{\gA}(a)$ and
	$\p_{\sigma^{-1}\circ\delta}^{\gA}(a')$ are atoms containing
	$s\circ\sigma$. Thus they are equal, and since
	$\p_{\sigma^{-1}\circ\delta}^{\gA}$ is an automorphism, $a=a'$.
	Consequently, the definition \eqref{eq:rep} assigns each $E(s,\sigma,\delta)$ in
	\eqref{eq:E-partition} to exactly one repetition-free atom.

	It follows immediately that the sets
	\begin{align}
		\{\rep(a):a\in\RfAt\}
		\label{eq:rep-partition}
	\end{align}
	are pairwise disjoint and their union is $F(S)$.
	They are also all nonempty. Indeed, let $a\in\RfAt$ and choose
	$z\in a$. Since $z\in S$, there are $s\in R$ and
	$\sigma\in PT_\alpha$ such that $z=s\circ\sigma$. Taking
	$\delta=\sigma$ in \eqref{eq:rep}, we have
	\[
		s\circ\sigma\in
		\p_{\sigma^{-1}\circ\sigma}^{\gA}(a)=a,
	\]
	and therefore $\P^V_{\sigma^{-1}}(s_\sigma)\subseteq\rep(a)$.
	The set $s_\sigma$ is nonempty by Lemma~\ref{huha-lemma}, so
	$\rep(a)\neq\emptyset$.

	We now consider the atoms which are not repetition-free. By definition,
	for every $a\in At\setminus\RfAt$ there are $k\neq l$ such that
	$a\subseteq\d_{kl}^{\gA}$, and hence $G(a)=F(a)$. Since the atoms of
	$\gA$ partition ${}^{\alpha}U$ and $S=\bigcup\{a:a\in\RfAt\}$,
	the atoms outside $\RfAt$ partition ${}^{\alpha}U\setminus S$.
	The map $F$ is a Boolean embedding on the full set algebra, so
	\begin{align}
		\{F(a):a\in At\setminus\RfAt\}
		\label{eq:nonrf-partition}
	\end{align}
	is a partition of
	\[
		F({}^{\alpha}U\setminus S)
		={}^{\alpha}V\setminus F(S).
	\]
	Each of these sets is nonempty because every atom $a$ is nonempty and
	$F(a)\neq\emptyset$.

	Combining \eqref{eq:rep-partition} and
	\eqref{eq:nonrf-partition}, and recalling the definition of $G$ on
	atoms, we conclude that $\{G(a):a\in At\}$
	is a partition of ${}^{\alpha}V$ into nonempty sets.

	It remains to verify the preservation of operations. Define $G(x)$ by
	\eqref{eq:G-extension}. Since the $G(a)$'s form a partition, we have
	\[
		G(0)=\emptyset,\qquad G(1)={}^{\alpha}V.
	\]
	For $x,y\in A$, an atom $a$ is below $x\cup y$ if and only if
	$a\leq x$ or $a\leq y$. Hence 
	\[
		G(x\cup y)=G(x)\cup G(y).
	\]
	Similarly, for every atom $a$, exactly one of $a\leq x$ and
	$a\leq -x$ holds. Since the sets $G(a)$ partition ${}^{\alpha}V$,
	\[
		G(-x)={}^{\alpha}V\setminus G(x).
	\]
	Thus $G$ is a Boolean homomorphism.

	Finally, $G$ is injective. If $x\neq 0$, atomicity gives an atom
	$a\leq x$. Since $G(a)\neq\emptyset$,
	\[
		\emptyset\neq G(a)\subseteq G(x),
	\]
	and therefore $G(x)\neq\emptyset$. Thus the kernel of $G$ is
	$\{0\}$, and a Boolean homomorphism with trivial kernel is injective.
	\end{proof}

	Using Lemma \ref{lem:P-re-jo} and \eqref{eq:pPad}, both for repetition-free
	and non repetition-free atoms $a$, $G(\p_{\sigma}^{\gA}a) = \P_{\sigma}(G(a))$, 
	and thus $G$ is a homomorphism with respect to $\p$ as well. 
	
	$G$ is a homomorphism for the cylindrifications: we need that for
	any atom $a$, $G(\C_i^U(a)) = \C_i^V(G(a))$. By Lemma \ref{lem:C-re-jo}, 
	$\C_i^V(G(a))$ is always $\C_i^V(F(a))$. Now, if $z\in b$ for an atom 
	$b\subseteq\C_i^U(a)$, then we have two cases. Either $b$ is not
	repetition-free, in which case $F(\{z\})\subseteq G(b)$; or $b$ is 
	repetition-free. 
	In this latter case there are $s\in R$ and $\tau\in PT_\alpha$ such that
	$z\circ\tau=s$. Since $F(\{s\})=\hat s=\bigcup_{\sigma\in PT_\alpha}s_\sigma$,
	and $z=s\circ\tau^{-1}$, we have
	\[
		F(\{z\})=\P_\tau^V(F(\{s\})) =\bigcup_{\sigma\in PT_\alpha}\P_\tau^V(s_\sigma).
	\]
	For each $\sigma\in PT_\alpha$, let $c_\sigma$ be the atom for which
	$s\circ\sigma\in \p_{\sigma^{-1}}^{\gA}\p_{\tau^{-1}}^{\gA}(c_\sigma)$.
	Taking $\delta=\tau^{-1}$ in the definition of $\rep$ gives
	$\P_\tau^V(s_\sigma)\subseteq\rep(c_\sigma)$.
	Moreover, $c_\sigma\subseteq\C_i^U(a)$. Indeed, since
	$s\circ\tau^{-1}=z\in\C_i^U(a)$,
	\[
	s\circ\sigma
	 \in \P_{\sigma^{-1}\circ\tau^{-1}}^U(\C_i^U(a))
	 =\p_{\sigma^{-1}}^{\gA}\p_{\tau^{-1}}^{\gA}(\C_i^U(a)),
	\]
	where the last equality follows from the fact that the abstract and concrete
	permutations agree on cylindrified elements. Hence, by atomicity,
	$c_\sigma\subseteq\C_i^U(a)$. Since $c_\sigma$ is repetition-free,
	$G(c_\sigma)=\rep(c_\sigma)$, and therefore 
	\[
	F(\{z\})
	 \subseteq \bigcup_{\sigma\in PT_\alpha}G(c_\sigma)
	 \subseteq G(\C_i^U(a)).
	\]	
	Putting the cases together, we end up at
	\begin{align}
		&G(\C_i^U(a)) = \bigcup\big\{ G(b):\; b\subseteq \C_{i}^U(a), b\in At \big\}
		\supset \bigcup\big\{F(z): z\in \C_i^U(a) \big\}  \\
		&= F(\C_{i}^U(a)) = \C_{i}^V(F(a))
		\overset{\text{\eqref{eq:uccso}}}{=} 
		\C_{i}^V(G(a)) .
	\end{align}
	
	Finally, $G(\d_{ij}^{\gA}) = F(\d_{ij}^{\gA}) = \d_{ij}'$. As $\S_{ij}^V$
	can be expressed using $\C_{i}^V$ and $\d_{ij}'$ (cf. \eqref{uj:eq:s=cd}), 
	it follows that $G(\S_{ij}^U(a)) = \S_{ij}^V(G(a))$, and the proof is complete.
\end{proof}

As a corollary to Theorem \ref{thm:2a}, we can fill in the status of the open line 
 from $cs$ to $csp$  in \cite[Figure 1]{AAPALI97ec}: it has to be a dashed line indicating finite axiomatizability.

\begin{theorem}\label{thm:2b} 
	$\RA_{\alpha}^{csp}$ is finitely axiomatizable over $\RA_{\alpha}^{cs}$
	(by the equations (F0)--(F9) axiomatizing $\FPA_{\alpha}$).
\end{theorem}	
\begin{proof}
	This follows from \eqref{uj:*} and $\RA_{\alpha}^{csp}\subseteq\FPA_{\alpha}$.  
\end{proof}

We now prove Theorem~\ref{main3-thm}. The proof is completely analogous to that of Theorem~\ref{main2-thm}. 
By $\FRPA_{\alpha}=\mathbf{SP}\RA_{\alpha}^{csp}$ and $\SRd_{cs}\FRPA_{\alpha}=\mathbf{SP}\RA_{\alpha}^{cs}$, Theorem~\ref{main3-thm} states that $\gA\in\mathbf{SP}\RA_{\alpha}^{csp}$ iff ($\Rd_{cs}\gA\in\mathbf{SP}\RA_{\alpha}^{cs}$  and equations \eqref{main10}--\eqref{main13} hold in $\gA$). By Theorem~\ref{thm:2b} we get that $\gA\in\RA_{\alpha}^{csp}$ iff ($\gA\in\RA_{\alpha}^{cs}$  and equations (F0)--(F9) hold in $\gA$).  Equations \eqref{main10}--\eqref{main13} in Theorem \ref{main3-thm} are exactly the equations in (F0)--(F9) in which the permutation operations occur. Theorem~\ref{main3-thm} has been proved.\\

We prove Theorem~\ref{main1-thm}. By Theorems~\ref{main2-thm},~\ref{main3-thm} and the last item (FS) in Lemma~\ref{pc-lemma} we get that $\FRPA_{\alpha}$ is axiomatized over $\RDF_{\alpha}$ by the equations (F0)--(F9) defining $\FPA_{\alpha}$.  \cite[Theorem 1(i)]{Sain-Tho} states that $\FPA_{\alpha}$ is term-definitionally equivalent to $\PA_{\alpha}$. The same term-definitions take $\FRPA_{\alpha}$ to $\RPA_{\alpha}$. Hence, $\RPA_{\alpha}$ is finitely axiomatized over $\RDF_{\alpha}$ by the equations defining $\PA_{\alpha}$. Examining the defining equations of $\PA$, in e.g., \cite[Def.5.4.1]{HMT}, one sees that among these equations \eqref{main1}--\eqref{main4} are the ones in which the substitution operations $\s_{\sigma}$ occur and the rest follow easily from (F0)--(F9). Theorem~\ref{main1-thm} has been proved.\\

\noindent We conclude by recalling a question from \cite{AAPALI97ec}, which remains open.

\begin{problem}
	Are permutations finitely axiomatizable over the Boolean operations together with cylindrifications? 
	I.e., does the class $\RA^{cp}_{\alpha}$ admit a finite axiomatization over $\RA^{c}_{\alpha}$?
\end{problem}

\paragraph{{\bf Acknowledgement.}}
We are grateful to the anonymous referee for the exceptional care and effort he devoted to reading our manuscript. The referee's detailed comments not only helped us to correct and clarify several points, but also suggested a number of useful simplifications. In particular, the referee proposed an alternative route through the main argument which avoids the use of union-completeness. We have chosen to retain our original proof, suitably revised, but we believe that the referee's alternative argument is of independent interest and may be useful in later investigations. We therefore include a condensed version of it in the Appendix.

\appendix
\section{Alternative proofs}
\label{app:without-atoms}

The anonymous referee suggested alternative proofs for Theorems~\ref{thm:1a}
and~\ref{thm:2a} which avoid atoms and union-completeness. We give these
arguments in a condensed version and specify the changes to Sections~\ref{sec:main1}
and~\ref{sec:povercs}. Canonicity is still used to obtain a complete
algebra, so that Theorem~\ref{rdsc-thm}\ref{rdsc1} supplies the
diagonal elements, but its set representation need not preserve
infinite suprema. After the diagonals have been obtained, completeness
is not used either. Throughout this appendix $3\leq\alpha<\omega$.

\begin{lemma}\label{app:reduction}
To prove Theorems~\ref{thm:1a} and~\ref{thm:2a}, respectively, it suffices
to prove the following implications for \emph{complete} BAOs $\gA$:
\begin{align}
 \gA\in\FSC_\alpha,\quad \Rd_c\gA\in\RA^c_\alpha
 &\quad\Longrightarrow\quad \gA\in\RA^{cs}_\alpha,
 \label{app:reduction-c}\\
 \gA\in\FPA_\alpha,\quad \Rd_{cs}\gA\in\RA^{cs}_\alpha
 &\quad\Longrightarrow\quad \gA\in\RA^{csp}_\alpha.
 \label{app:reduction-cs}
\end{align}
\end{lemma}
\begin{proof}
	The proof is essentially the same as that of Lemma \ref{lem:1}. Canonicity
	is used to get a complete extension.
\end{proof}

Accordingly, in Lemma~\ref{lem:1} replace \eqref{*} by
\eqref{app:reduction-c}, quantified over complete BAOs, and
in the proof of Theorem~\ref{thm:2a} replace \eqref{uj:*} by
\eqref{app:reduction-cs} with the same qualification. 
At the subsequent passage to a set
algebra, take any set representation of the indicated reduct and
transport the remaining operations along it. Completeness as an
abstract BAO is preserved by this isomorphism; but union-completeness 
is not asserted. In Section~\ref{sec:povercs},
start likewise with a complete $\gA$ satisfying 
\eqref{app:reduction-cs} and represent its $cs$-reduct.
Theorem~\ref{rdsc-thm}\ref{rdsc1} now gives exactly the diagonal
expansions and equations \eqref{eq:s=cd} and \eqref{uj:eq:s=cd} used
in the two sections. 

For the next theorem no completeness or atomicity is assumed. 

\begin{theorem}\label{app:diagonals}
Suppose that $\gA=\langle A,\cup,\setminus,\C_i^U,
\s_{ij}^{\gA}\rangle_{i,j<\alpha}$ has a cylindric expansion by
elements $\d_{ij}^{\gA}\in A$ satisfying \eqref{eq:s=cd}. For every $i,j<\alpha$,
\begin{equation}
 \D_{ij}^U\subseteq\d_{ij}^{\gA}
 \quad\Longleftrightarrow\quad
 \S_{ij}^U(x)=\s_{ij}^{\gA}(x)\text{ for every }x\in A.
 \label{app:diagonal-criterion}
\end{equation}
Moreover, $\gA\in\RA^{cs}_\alpha$.
\end{theorem}
\begin{proof}
For $i=j$, both sides of \eqref{app:diagonal-criterion} hold. If
$i\neq j$ and $\D_{ij}^U\subseteq\d_{ij}^{\gA}$, then
$\S_{ij}^U(x)=\C_i^U(\D_{ij}^U\cap x)
\subseteq\C_i^U(\d_{ij}^{\gA}\cap x)=\s_{ij}^{\gA}(x)$.
Apply this also to $-x$. Since both substitution operations preserve
complements, the reverse inclusion follows. Conversely, equality of
the substitution operations gives
$\S_{ij}^U(-\d_{ij}^{\gA})
=\s_{ij}^{\gA}(-\d_{ij}^{\gA})=\emptyset$.
Every member of $\D_{ij}^U\setminus\d_{ij}^{\gA}$ would belong to
this empty set, proving the other direction.

We construct a representation satisfying the left side of
\eqref{app:diagonal-criterion}. The case of the trivial 
algebra is straightforward, so we 
assume $U\neq\emptyset$. Retain the enlargement
$V=U\times W$, where $W$ is infinite and $|W|>|U|$, and the embedding
$F$ defined using \eqref{eq:hats}. Thus $\d_{ij}'=F(\d_{ij}^{\gA})$,
and the algebra $\gA'$ has concrete cylindrifications.
For $0<i<\alpha$ consider the relation
$\dom_{0i}(\d_{0i}')=
\{\langle s(0),s(i)\rangle:s\in\d_{0i}'\}$ on $V$.
Lemma~\ref{egeszV}, whose proof uses only cylindric identities, gives
both its domain and range as $V$. Each element has $|W|=|V|$ possible
neighbours in either direction. Indeed, once the first coordinates of
a pair come from a sequence in $\d_{0i}^{\gA}$, the two $W$-coordinates
can be chosen independently in its image under $F$.

This can be done by a back-and-forth argument: Let $X$ be an
infinite set of cardinality $\kappa$ and assume $E\subseteq X\times X$ gives
each element $\kappa$ neighbours in each direction, that is, for each $b,c\in X$,
\[
	|\{x\in X: (b,x)\in E \}| = |\{x\in X: (x,c)\in E \}| = |X| = \kappa\,.
\]
We show that $E$ contains a bijection $\beta:X\to X$.
Enumerate
$X=\{x_\xi:\xi<\kappa\}$. Starting with the empty partial bijection,
at stage $\xi$ first put $x_\xi$ into the domain, if necessary, by
choosing an unused $E$-neighbour as its image; then put $x_\xi$ into the
range, if necessary, by choosing an unused $E$-neighbour as its
preimage. Take unions at limit stages. Before each choice, fewer than
$\kappa$ elements have been used, so both choices are possible.
The union is a bijection $X\to X$ contained in $E$. 

Apply this to $E=\dom_{0i}(\d_{0i}')$ to obtain a bijection
$f_{0i}:V\to V$ with $f_{0i}\subseteq\dom_{0i}(\d_{0i}')$;
put $f_{00}=\id$. Define $f(s)(i)=f_{0i}^{-1}(s(i))$ and
$G(X)=f[X]$, as in Lemma~\ref{lemma:n0}. The map $G$ preserves the
Boolean operations and the concrete cylindrifications: changing the
$i$th coordinate before applying $f$ allows exactly the same values
as changing it afterwards, since $f_{0i}$ is a bijection.
Write $\d_{ij}''=G(\d_{ij}')$.

If $s(0)=s(i)$ and $z=f^{-1}(s)$, then
$z(i)=f_{0i}(z(0))$. The inclusion
$f_{0i}\subseteq\dom_{0i}(\d_{0i}')$, together with
$\C_k^V\d_{0i}'=\d_{0i}'$ for $k\notin\{0,i\}$, implies
$z\in\d_{0i}'$. Hence $s\in\d_{0i}''$, proving
$\D_{0i}^V\subseteq\d_{0i}''$. Symmetry gives the case $(i,0)$;
for distinct nonzero $i,j$, the cylindric axioms give
\begin{equation*}
 \D_{ij}^V=\C_0^V(\D_{i0}^V\cap\D_{0j}^V)
 \subseteq\C_0^V(\d_{i0}''\cap\d_{0j}'')=\d_{ij}''.
\end{equation*}
The case $i=j$ is immediate. The transported substitutions still
satisfy \eqref{eq:s=cd2}, so \eqref{app:diagonal-criterion} shows
that they are $\S_{ij}^V$. Thus $G\circ F$ represents $\gA$.
\end{proof}

Now, these are the changes in Section~\ref{sec:main1} to get the
alternative proof: Keep the
construction through \eqref{eq:hats} and Lemma~\ref{egeszV}, but remove
the declaration of $At$, Lemma~\ref{dom-particio}, and the construction of the maps $f_{0i}^a$. Replace that construction and Lemma~\ref{lemma:n1}
by the back-and-forth construction above. Keep the definitions
of $f$ and $G$ and Lemmas~\ref{lemma:n0}--\ref{lemma:n2}, taking
$G(X)=f[X]$ and using
$f_{0i}\subseteq\dom_{0i}(\d_{0i}')$ in place of
\eqref{eq:star}. The proof of Lemma~\ref{lemma:n3} simplifies to
\eqref{app:diagonal-criterion} applied to $\gA''$.
The assertions of atomicity and union-completeness following the
constructions of $F$ and $G$ are then unnecessary. 
Together with Lemma~\ref{app:reduction}, this proves
Theorems~\ref{thm:1a} and~\ref{thm:2}.\\

As for Section~\ref{sec:povercs}, keep everything intact until 
Lemma \ref{uj:pij-igazi-dij-alatt}. The following single 
theorem gives the rest of the construction after Lemma \ref{uj:pij-igazi-dij-alatt}. 
The theorem does not 
make use of completeness or atomicity.

\begin{theorem}\label{app:permutations}
Suppose that
$\gA=\langle A,\cup,\setminus,\C_i^U,\S_{ij}^U,
\p_{ij}^{\gA}\rangle_{i,j<\alpha}\in\FPA_\alpha$ has a cylindric
expansion by elements $\d_{ij}^{\gA}\in A$ satisfying
\eqref{uj:eq:s=cd}. Then $\gA\in\RA^{csp}_\alpha$.
\end{theorem}
\begin{proof}
The case $U=\emptyset$ is immediate, so assume $U\neq\emptyset$.
Use the Boolean algebra $\gB$ from
Lemma~\ref{uj:pij-igazi-Bn}. It is invariant under every
$\p_\sigma^{\gA}$: by \ref{lemma1.6-kov}, each generator
$\C_i^U(x)$ is sent to the generator
$\C_{\sigma(i)}^U(\p_\sigma^{\gA}x)$.
Consequently, Lemma~\ref{uj:pij-igazi-Bn} extends from transpositions
to every permutation, giving \eqref{eq:aha} for all $b\in B$.
Moreover, $\d_{kl}^{\gA}\in B$, since it is closed under any
cylindrifications $\C_i$ with $i\notin\{k,l\}$. Every $a\in A$ below
$\d_{kl}^{\gA}$, with $k\neq l$, belongs to $B$ as well:
$a=\C_k^U(a)\cap\d_{kl}^{\gA}$, as proved in
Lemma~\ref{uj:pij-igazi-dij-alatt}. In particular, if
$\d^{\gA}=\bigcup_{k\neq l<\alpha}\d_{kl}^{\gA}$ and
$e^{\gA}=-\d^{\gA}$, then both are in $B$, and every element of $A$
below $\d^{\gA}$ is in $B$, by distributing over this finite union.

Put $W=PT_\alpha$ and $V=U\times W$, and let $F$ as in \eqref{def:F}. 
Give $W$ any group operation $\star$ as in Lemma~\ref{huha-lemma}. For $\tau\in PT_\alpha$
and $X\subseteq{}^\alpha U$, put
\begin{equation}
 F_\tau(X)=\{s\times w:s\in X,\ w\in{}^\alpha W,
                \ w_0\star\cdots\star w_{\alpha-1}=\tau\}.
 \label{app:Ftau}
\end{equation}
Thus $F_\tau(\{s\})=s_\tau$ in the notation of
Lemma~\ref{huha-lemma}, and $F(X)=\bigcup_\tau F_\tau(X)$.
For $\sigma\in PT_\alpha$, put
$Q_\sigma(a)=\P_\sigma^U(\p_{\sigma^{-1}}^{\gA}(a))$ for $a\in A$.
Each $Q_\sigma:A\to\Sb({}^\alpha U)$ preserves Boolean operations, but 
note that its values need not belong to $A$.
Choose $R\subseteq{}^\alpha U$ containing exactly one member of
each set $PT(s)=\{s\circ\sigma:\sigma\in PT_\alpha\}$ for
every $s\in{}^\alpha U$. Define, for every $a\in A$,
\begin{equation}
 \rep(a)=\bigcup_{\sigma,\tau\in PT_\alpha}
 \P_{\sigma^{-1}}^V
 F_\tau\bigl(R\cap Q_\tau(\p_\sigma^{\gA}(a))\bigr)\quad\in \Sb({}^{\alpha}V).
 \label{app:rep}
\end{equation}
Now, for $X\subseteq{}^\alpha U$, $a\in A$, and $i<\alpha$, we have
\begin{align}
 \C_i^V F_\tau(X)&=F(\C_i^U X)=\C_i^V F(X),
 \label{app:F-cyl}\\
 \C_i^U Q_\sigma(a)&=Q_\sigma(\C_i^U a),
 \label{app:Q-cyl}\\
 \P_{\sigma^{-1}}^U Q_\tau(\p_\sigma^{\gA}(a))
 &=Q_{\sigma^{-1}\tau}(a).
 \label{app:Q-compose}
\end{align}
For the first equality, once all other $W$-coordinates have been
fixed, there is exactly one value of the $i$th coordinate making
their ordered product equal to $\tau$. If the products before and
after that coordinate are $u$ and $v$, this value is
$u^{-1}\star\tau\star v^{-1}$, where inverses and empty products
are taken in $(W,\star)$. This proves
\eqref{app:F-cyl}.
Equation~\eqref{app:Q-cyl} follows by applying
\ref{lemma1.6-kov} and the concrete identity
$\C_i^U\P_\sigma^U=\P_\sigma^U\C_{\sigma^{-1}(i)}^U$.
Equation~\eqref{app:Q-compose} follows from the composition laws
for $\P$ and $\p$. We will also use that $F$ commutes with the
concrete permutations and preserves set-theoretic unions.

For $b\in B$, equation~\eqref{eq:aha} gives $Q_\xi(b)=b$ for all
$\xi$. Hence \eqref{app:Q-compose} implies
$Q_\tau(\p_\sigma^{\gA}b)=\P_\sigma^U b$, and therefore
\begin{align}
 \rep(b)
 &=\bigcup_\sigma\P_{\sigma^{-1}}^V
                  F(R\cap\P_\sigma^U b)\notag\\
 &=F\left(\left(\bigcup_\sigma\P_{\sigma^{-1}}^U R\right)
                  \cap b\right)=F(b),\qquad b\in B.
 \label{app:rep-on-B}
\end{align}
Here $\bigcup_\sigma\P_{\sigma^{-1}}^U R={}^\alpha U$ by the
choice of $R$. In particular, \eqref{app:rep-on-B} applies to
$\C_i^U(a)$, $\d_{kl}^{\gA}$, $\d^{\gA}$, and $e^{\gA}$.

We first verify preservation of cylindrifications. Using
\eqref{app:F-cyl}--\eqref{app:Q-compose} gives
\begin{align*}
 \C_i^V\rep(a)
 &=\bigcup_{\sigma,\tau}\P_{\sigma^{-1}}^V
    F\bigl(\C_{\sigma(i)}^U
       (R\cap Q_\tau(\p_\sigma^{\gA}a))\bigr)\\
 &=F\left(\C_i^U\bigcup_{\sigma,\tau}
       (\P_{\sigma^{-1}}^U R\cap Q_{\sigma^{-1}\tau}(a))\right)\\
 &=F\left(\C_i^U\bigcup_\xi Q_\xi(a)\right)
  =F\left(\bigcup_\xi Q_\xi(\C_i^U a)\right)\\
 &=F(\C_i^U a)=\rep(\C_i^U a).
\end{align*}
In the third equality, for each fixed $\sigma$ the index
$\xi=\sigma^{-1}\tau$ ranges over all permutations, and the
sets $\P_{\sigma^{-1}}^U R$ cover ${}^\alpha U$. The last line
uses $\C_i^U a\in B$ and \eqref{app:rep-on-B}.

The map $\rep$ preserves finite unions and $0$ by definition, and
$\rep(1)=F(1)={}^\alpha V$ by \eqref{app:rep-on-B}.
We prove that disjoint elements have disjoint images.
On elements below $\d^{\gA}$, the map is $F$; also
$\rep(\d^{\gA})=F(\d^{\gA})$ and
$\rep(e^{\gA})=F(e^{\gA})$ are disjoint. Splitting each element
$a$ as $(a\cap\d^{\gA})\cup(a\cap e^{\gA})$, it therefore
suffices to consider disjoint $a,b\in A$ below $e^{\gA}$.
Suppose $t\times w\in\rep(a)\cap\rep(b)$. Monotonicity and
\eqref{app:rep-on-B} give $t\in e^{\gA}$, so $t$ is
repetition-free by \eqref{eq:Dij:dij}.
There is consequently a unique $\sigma$ such that
$t\circ\sigma^{-1}\in R$: the representative is unique, and
the coordinates of $t$ are distinct. In any occurrence of
$t\times w$ in the union \eqref{app:rep}, this must be the
index $\sigma$, and $\tau$ must be the product of the
coordinates of $w\circ\sigma^{-1}$, by \eqref{app:Ftau}.
Membership in both images therefore implies
\begin{equation*}
 t\circ\sigma^{-1}\in
 Q_\tau(\p_\sigma^{\gA}a)\cap Q_\tau(\p_\sigma^{\gA}b)
 =Q_\tau(\p_\sigma^{\gA}(a\cap b))=\emptyset,
\end{equation*}
a contradiction. This proves disjointness preservation, and hence
preservation of complements.

To prove injectivity, take nonempty $a\in A$ and choose
$\sigma$ with $R\cap\P_\sigma^U a\neq\emptyset$.
Since $Q_\sigma\p_\sigma^{\gA}=\P_\sigma^U$, the term with
$\tau=\sigma$ in \eqref{app:rep} is
$\P_{\sigma^{-1}}^V F_\sigma(R\cap\P_\sigma^U a)$.
It is nonempty by the group equation used in
\eqref{app:F-cyl}, so $\rep(a)\neq\emptyset$. Thus this
Boolean homomorphism has trivial kernel.

For any permutation $\theta$, reindexing by
$\xi=\sigma\theta^{-1}$ yields
\begin{align*}
 \P_\theta^V\rep(a)
 &=\bigcup_{\sigma,\tau}\P_{\theta\sigma^{-1}}^V
      F_\tau(R\cap Q_\tau(\p_\sigma^{\gA}a))\\
 &=\bigcup_{\xi,\tau}\P_{\xi^{-1}}^V
      F_\tau(R\cap Q_\tau(\p_\xi^{\gA}\p_\theta^{\gA}a))
  =\rep(\p_\theta^{\gA}a).
\end{align*}
Finally, $\rep(\d_{ij}^{\gA})=F(\d_{ij}^{\gA})$, and
$\D_{ij}^V\subseteq F(\D_{ij}^U)
\subseteq F(\d_{ij}^{\gA})$. Since $\rep$ preserves the
Boolean operations and cylindrifications, the transported
substitutions satisfy \eqref{uj:eq:s=cd} with
$\rep(\d_{ij}^{\gA})$ in place of $\d_{ij}^{\gA}$.
Apply \eqref{app:diagonal-criterion} to the transported cylindric
expansion: these substitutions are exactly $\S_{ij}^V$.
Thus $\rep$ is the required representation of $\gA$.
\end{proof}

\end{document}